# A New Approach in Solving Regular and Singular Conformable Fractional Coupled Burger's Equations

AMJAD E. HAMZA[1], ABDELILAH K. SEDEEG[2,3], RANIA SAADEH[4]*, AHMAD QAZZA[4], RAED KHALIL[5]

[1]Department of Mathematic, Faculty of Sciences,
University of Ha'il, Ha'il 2440,
SAUDI ARABIA

[2]Department of Mathematics, Faculty of Education,
Holy Quran and Islamic Sciences University, Sudan,
SUDAN

[3]Department of Mathematics, Faculty of Sciences and Arts-Almikwah,
Albaha University,
SAUDI ARABIA

[4]Department of Mathematics, Faculty of Science,
Zarqa University, Zarqa 13110,
JORDAN

[5]Department of Computer Information Systems, Faculty of Prince Abdullah Bin Ghazi,
Balqaa Applied Univesity,
JORDAN

*Corresponding Author*

*Abstract:* - The conformable double ARA decomposition approach is presented in this current study to solve one-dimensional regular and singular conformable functional Burger's equations. We investigate the conformable double ARA transform's definition, existence requirements, and some basic properties. In this study, we introduce a novel interesting method that combines the double ARA transform with Adomian's decomposition method, in order to find the precise solutions of some nonlinear fractional problems. Moreover, we use the new approach to solve Burgers' equations for both regular and singular conformable fractional coupled systems. We also provide several instances to demonstrate the usefulness of the current study. Mathematica software has been used to get numerical results.

## 1 Introduction

Fractional partial differential equations have drawn significant interest from a wide range of specialists in applied sciences and engineering, including acoustics, control, and viscoelasticity. In many areas of mathematics and physics, partial differential equations are crucial. To examine several time-fractional partial differential equations, the authors use a novel strategy termed the "simplest equation method", [1], [2].

In the context of applied sciences like mathematical modeling and fluid mechanics, this work focuses on Burger's equation. Burger's equation was initially brought up about steady-state solutions, in fact, [3]. Burger later changed the approach to characterize the viscosity of certain fluid types, [4]. The conformable double Laplace transform approach, which was first proposed in [5], was improved and used to solve fractional partial differential equations. This method has been used by a number of academics to obtain precise and numerical solutions to this type of equation. To precisely solve time-fractional Burger's equations, other scholars used the first integral technique, [6]. Another set of researchers developed the coupled Burger's equation solution using the generalized two-dimensional differential transform approach, [7], [8], [9], [10].

The ARA transform is a revolutionary integral transform that Saadeh and others introduced in, [11]. ARA is a new transform; it is not an acronym. It has novel properties, including the ability to generate numerous transforms by varying the value of the index n, a duality with the Laplace transform, and the capacity to get around the singularity at time zero. Many researchers have studied the new approach and implemented it to solve many problems by merging it with other numerical methods or other transforms, such as ARA-Sumudu transform, [12], [13], Laplace-ARA transform, [14], double ARA transform, [15], [16], ARA residual power series method, [17], [18].

In this article, we choose to build a unique combination of Adomian's decomposition method and the double ARA transform, so that we obtain the advantages of these two methods and fully utilize these two potent techniques. The conformable double ARA transform method will be introduced in this research in combination with Adomian's decomposition method, [19], to solve systems of conformable fractional partial differential equations.

With the help of the conformable double ARA decomposition approach, this study aims to provide analytical solutions for the coupled, one-dimensional, singular, and regular conformable fractional Burger's equations (CDARADM). The following space-time fractional order coupled with Burger's equations were described in [20], and are given below:

$$\begin{aligned}\frac{\partial^q u}{\partial t^q}-\frac{\partial^{2p} u}{\partial x^{2p}}+\lambda u\frac{\partial^p u}{\partial x^p}+\alpha\frac{\partial^p}{\partial x^p}(uv)&=k\left(\frac{x^p}{p},\frac{t^q}{q}\right),\\ \frac{\partial^q v}{\partial t^q}-\frac{\partial^{2p} v}{\partial x^{2p}}+\lambda v\frac{\partial^p v}{\partial x^p}+\beta\frac{\partial^p}{\partial x^p}(uv)&=l\left(\frac{x^p}{p},\frac{t^q}{q}\right).\end{aligned} \tag{1}$$

This article is organized as follows, in the following section, we present the ARA transform with the main characteristics. In Section 3, we introduce some preliminaries about the conformable fractional derivatives. The conformable ARA transform and some related results are presented in Section 4. In Section 5, we introduce some numerical experiments to prove the efficiency and applicability of the new method.

# 2 ARA Integral Transforms, [11]

**Definition 1.** If $h(x)$ is a continuous function on $(0,\infty)$, then the ARA transform of order $n$

$$\begin{aligned}\mathcal{G}_n[h(x)](r)&=Q(n,r)\\&=r\int_0^\infty x^{n-1}e^{-rx}h(x)\,dx,\quad r>0,\end{aligned} \tag{2}$$

and the inverse ARA transform is defined as

$$\begin{aligned}&\mathcal{G}_{n+1}^{-1}\left[\mathcal{G}_{n+1}[h(x)]\right]\\&=\frac{(-1)^{2n}}{2\pi i}\int_{c-i\infty}^{c+i\infty}e^{rx}\,Q(r)dr=h(x),\end{aligned} \tag{3}$$

where

$$Q(r)=\int_0^\infty e^{-rx}h(x)\,dx.$$

**Theorem 1.** Let $h(x)$ be piecewise continuous in every finite interval $0\le x\le\alpha$ and satisfies the condition

$$|x^{n-1}h(x)|\le Me^{\alpha x}.$$

where $M$ is a positive constant. Then, ARA integral transform exists for all $r>\alpha$.

**Proof.** Using the definition of ARA transform, we get

$$|Q(n,r)|=\left|r\int_0^\infty x^{n-1}e^{-rx}h(x)\,dx\right|.$$

Thus, we get

$$\begin{aligned}|Q(n,r)|&=\left|r\int_0^\beta x^{n-1}e^{-rx}h(x)\,dx\right.\\&\qquad\left.+r\int_\beta^\infty x^{n-1}e^{-rx}h(x)\,dx\right|\\&\le r\left|\int_\beta^\infty x^{n-1}e^{-rx}h(x)\,dx\right|\\&\le r\int_\beta^\infty e^{-rx}|x^{n-1}h(x)|\,dx\\&\le r\int_\beta^\infty e^{-rx}Me^{\alpha x}\,dx\\&=rM\int_\beta^\infty e^{-(r-\alpha)x}dx\\&=\frac{rM}{r-\alpha}e^{-\beta(r-\alpha)}.\end{aligned}$$

Hence, the integral exists for all $r>\alpha$, and $\mathcal{G}_{n+1}[h(x)]$ exists.

Now, we present some properties of ARA transform.

If $H(n,r) = \mathcal{G}_n[h(x)]$ and $G(n,r) = \mathcal{G}_n[g(x)]$ and $a, b \in \mathbb{R}$, then

- $$\mathcal{G}_n[a\,h(x) + b\,g(x)] = a\,\mathcal{G}_n[h(x)] + b\,\mathcal{G}_n[g(x)]. \quad (4)$$
- $$\mathcal{G}_n{}^{-1}[a\,H(n,r) + b\,G(n,r)] = a\,\mathcal{G}_n{}^{-1}[H(n,r)] + b\,\mathcal{G}_n{}^{-1}[G(n,r)]. \quad (5)$$
- $$\mathcal{G}_n[x^\alpha] = \frac{\Gamma(\alpha+n)}{r^{\alpha+n-1}}, \quad \alpha > 0. \quad (6)$$
- $$\mathcal{G}_n[e^{ax}] = \frac{r\Gamma(n)}{(r-a)^n}, \quad a \in \mathbb{R}. \quad (7)$$
- $$\mathcal{G}_n[\sin ax] = \frac{r}{2i}\Gamma(n)\left(\frac{1}{(r-ia)^n} - \frac{1}{(r+ia)^n}\right), a \in \mathbb{R}. \quad (8)$$
- $$\mathcal{G}_n\left[h^{(n)}(x)\right] = (-1)^{n-1} r \frac{d^{n-1}}{ds^{n-1}}\left(r^{n-1}\mathcal{G}_1[h(x)] - \sum_{k=1}^{n} r^{n-k}\, h^{(k-1)}(0)\right), \quad (9)$$

where $\mathcal{G}_1[h(x)]$ is the ARA transform of a continuous function $h(x)$ of order one on $[0,\infty)$, and it is given by

$$\mathcal{G}_1[h(x)] = H(r) = r\int_0^\infty e^{-rx} h(x)\,dx, \quad r > 0. \quad (10)$$

In this research, we denote $\mathcal{G}_1[h(x)]$ for $\mathcal{G}[h(x)]$.

# 3 Conformable Fractional Derivatives (CFD), [21], [22], [23], [24]

Conformable fractional derivatives are investigated and expanded in [15], [16], respectively. The following definitions of CFD are utilized in this study.

**Definition 2.** The CFD of $h\left(\frac{x^p}{p}\right)$ order $p$, where $h: (0,\infty) \to \mathbb{R}.$, is given by

$$\frac{d^p}{dx^p} h\left(\frac{x^p}{p}\right) = \lim_{\delta\to 0} \frac{h\left(\frac{x^p}{p} + \delta x^{1-p}\right) - h\left(\frac{x^p}{p}\right)}{\delta}, \quad (11)$$

$$\frac{x^p}{p} > 0, \quad 0 < p \le 1.$$

**Definition 3.** The conformable space fractional partial derivative of $h\left(\frac{x^p}{p}, \frac{t^q}{q}\right)$ order $p$, where $h\left(\frac{x^p}{p}, \frac{t^q}{q}\right): \mathbb{R} \times (0,\infty) \to \mathbb{R}$, is defined as

$$\frac{\partial^p}{\partial x^p} h\left(\frac{x^p}{p}, \frac{t^q}{q}\right) = \lim_{\delta\to 0} \frac{h\left(\frac{x^p}{p} + \delta x^{1-p}, \frac{t^q}{q}\right) - h\left(\frac{x^p}{p}, \frac{t^q}{q}\right)}{\delta}, \quad (12)$$

$$\frac{x^p}{p}, \frac{t^q}{q} > 0, \quad 0 < p, q \le 1.$$

**Definition 4.** Let $h\left(\frac{x^p}{p}, \frac{t^q}{q}\right): \mathbb{R} \times (0,\infty) \to \mathbb{R}$. Then, the conformable fractional partial derivative of $h\left(\frac{x^p}{p}, \frac{t^q}{q}\right)$ of order, $q$ is given by

$$\frac{\partial^q}{\partial t^q} h\left(\frac{x^p}{p}, \frac{t^q}{q}\right) = \lim_{\varepsilon\to 0} \frac{h\left(\frac{x^p}{p}, \frac{t^q}{q} + \xi t^{1-q}\right) - h\left(\frac{x^p}{p}, \frac{t^q}{q}\right)}{\xi}, \quad (13)$$

$$\frac{x^p}{p}, \frac{t^q}{q} > 0, \quad 0 < p, q \le 1.$$

## 3.1 Conformable Fractional Derivatives of Some Basic Functions

In the following arguments, we introduce the CFD for some basic functions.

i. Let $h\left(\frac{x^p}{p}, \frac{t^q}{q}\right) = \frac{x^p}{p}\frac{t^q}{q}$ . Then
   - $\frac{\partial^p}{\partial x^p}\left(\frac{x^p}{p}\frac{t^q}{q}\right) = \frac{t^q}{q}$,
   - $\frac{\partial^q}{\partial t^q}\left(\frac{x^p}{p}\frac{t^q}{q}\right) = \frac{x^p}{p}$.

ii. Let $h\left(\frac{x^p}{p}, \frac{t^q}{q}\right) = \left(\frac{x^p}{p}\right)^n$. Then
   - $\frac{\partial^p}{\partial x^p}\left(\frac{x^p}{p}\right)^n = n\left(\frac{x^p}{p}\right)^{n-1}$,
   - $\frac{\partial^q}{\partial t^q}\left(\frac{x^p}{p}\right)^n = 0.$

iii. Let $h\left(\frac{x^p}{p}, \frac{t^q}{q}\right) = \left(\frac{x^p}{p}\right)^n \left(\frac{t^q}{q}\right)^m$ . Then
   - $\frac{\partial^p}{\partial x^p}\left(\frac{x^p}{p}\right)^n\left(\frac{t^q}{q}\right)^m = n\left(\frac{x^p}{p}\right)^{n-1}\left(\frac{t^q}{q}\right)^m$,
   - $\frac{\partial^q}{\partial t^q}\left(\frac{x^p}{p}\right)^n\left(\frac{t^q}{q}\right)^m = m\left(\frac{x^p}{p}\right)^n\left(\frac{t^q}{q}\right)^{m-1}$.

iv. Let $h\left(\frac{x^p}{p}, \frac{t^q}{q}\right) = \sin\lambda\left(\frac{x^p}{p}\right)\sin\beta\left(\frac{t^q}{q}\right)$. Then
   - $\frac{\partial^p}{\partial x^p}\left(\sin\lambda\left(\frac{x^p}{p}\right)\sin\beta\left(\frac{t^q}{q}\right)\right) = \lambda\cos\lambda\left(\frac{x^p}{p}\right)\sin\beta\left(\frac{t^q}{q}\right)$,
   - $\frac{\partial^q}{\partial t^q}\left(\sin\lambda\left(\frac{x^p}{p}\right)\sin\beta\left(\frac{t^q}{q}\right)\right) = \beta\sin\lambda\left(\frac{x^p}{p}\right)\cos\beta\left(\frac{t^q}{q}\right)$.

v. Let $h\left(\frac{x^p}{p},\frac{t^q}{q}\right)=e^{\lambda\frac{x^p}{p}+\beta\frac{t^q}{q}}$. Then

- $\frac{\partial^p}{\partial x^p}\left(e^{\lambda\frac{x^p}{p}+\beta\frac{t^q}{q}}\right)=\lambda e^{\lambda\frac{x^p}{p}+\beta\frac{t^q}{q}}$,
- $\frac{\partial^q}{\partial t^q}\left(e^{\lambda\frac{x^p}{p}+\beta\frac{t^q}{q}}\right)=\beta e^{\lambda\frac{x^p}{p}+\beta\frac{t^q}{q}}$.

**Property 1**. If $h(x,t)$ is a differentiable function of order $\alpha$ and $\beta$ at the points $x$ and $t>0$, where $0<\alpha,\beta\leq 1$. Then

$$\frac{\partial^p}{\partial x^p}h\left(\frac{x^p}{p},\frac{t^q}{q}\right)= x^{1-p}\,\frac{\partial}{\partial x}h\left(\frac{x^p}{p},\frac{t^q}{q}\right),$$

$$\frac{\partial^q}{\partial t^q}h\left(\frac{x^p}{p},\frac{t^q}{q}\right)= t^{1-q}\,\frac{\partial}{\partial t}h\left(\frac{x^p}{p},\frac{t^q}{q}\right).$$

**Proof**. Using Definition **3** and putting $k=\delta\,x^{1-p}$ in Equation (12), we get

$$\begin{aligned}\frac{\partial^p}{\partial x^p}h\left(\frac{x^p}{p},\frac{t^q}{q}\right)&=\lim_{\delta\to 0}\frac{h\left(\frac{x^p}{p}+\delta x^{1-p},\frac{t^q}{q}\right)-\left(\frac{x^p}{p},\frac{t^q}{q}\right)}{\delta}\\ &=\lim_{k\to 0}\frac{h\left(\frac{x^p}{p}+k,\frac{t^q}{q}\right)-\left(\frac{x^p}{p},\frac{t^q}{q}\right)}{kx^{p-1}}\\ &= x^{1-p}\lim_{k\to 0}\frac{h\left(\frac{x^p}{p}+k,\frac{t^q}{q}\right)-\left(\frac{x^p}{p},\frac{t^q}{q}\right)}{k}\\ &=x^{1-p}\frac{\partial}{\partial x}h\left(\frac{x^p}{p},\frac{t^q}{q}\right).\end{aligned}$$

Similarly, we can easily prove that

$$\frac{\partial^q}{\partial t^q}h\left(\frac{x^p}{p},\frac{t^q}{q}\right)= t^{1-q}\,\frac{\partial}{\partial t}h\left(\frac{x^p}{p},\frac{t^q}{q}\right).$$

# 4 Conformable Double ARA Transform (CDARAT)

In this part of this study, we present the conformable double ARA transform using the following definitions.

**Definition 5.** Assume that $h(x)$ is a real-valued function defined on $[0,\infty)$ to $\mathbb{R}$, then the conformable ARA transform of $h\left(\frac{x^p}{p}\right)$ is given by

$$\mathcal{G}_x^p\left[h\left(\frac{x^p}{p}\right)\right] =r\int_0^\infty e^{-r\frac{x^p}{p}}\left[h\left(\frac{x^p}{p}\right)\right]x^{p-1}dx, \qquad (14)$$

$\forall r>0$,

If the integral exists.

**Definition 6.** Let $h\left(\frac{x^p}{p},\frac{t^q}{q}\right)$ be a piecewise continuous function on the interval $[0,\infty)\times[0,\infty)$ of exponential order, that is considered for some $a,b\in\mathbb{R}$, $\frac{x^p}{p},\frac{t^q}{q}>0$, $\frac{\left|h\left(\frac{x^p}{p},\frac{t^q}{q}\right)\right|}{e^{a\frac{x^p}{p}+b\frac{t^q}{q}}}$.Under these conditions the conformable double ARA transform is given by

$$\begin{aligned}\mathcal{G}_x^p\mathcal{G}_t^q\left[h\left(\frac{x^p}{p},\frac{t^q}{q}\right)\right]&=G_{p.q}(r,s)\\ &=rs\int_0^\infty\int_0^\infty e^{-r\frac{x^p}{p}-s\frac{t^q}{q}}h\left(\frac{x^p}{p},\frac{t^q}{q}\right)x^{p-1}t^{q-1}dx\,dt,\end{aligned} \qquad (15)$$

where $r,s\in\mathbb{C}$, $0<p,q\leq 1$ and the integrals with respect to $x^p$ and $t^q$ respectively, are taken by conformable fractional derivative.

**Property 2.** Let $h\left(\frac{x^p}{p},\frac{t^q}{q}\right)=f\left(\frac{x^p}{p}\right)g\left(\frac{t^q}{q}\right), x>0, t>0$. Then

$$\mathcal{G}_x^p\mathcal{G}_t^q\left[h\left(\frac{x^p}{p},\frac{t^q}{q}\right)\right]=\mathcal{G}_x[f(x)]\,\mathcal{G}_t[g(t)]. \qquad (16)$$

**Proof.** The definition of CDARAT implies

$$\begin{aligned}&\mathcal{G}_x^p\mathcal{G}_t^q\left[f\left(\frac{x^p}{p}\right)g\left(\frac{t^q}{q}\right)\right]\\ &=rs\int_0^\infty\int_0^\infty e^{-r\frac{x^p}{p}-s\frac{t^q}{q}}f\left(\frac{x^p}{p}\right)g\left(\frac{t^q}{q}\right)x^{p-1}t^{q-1}dx\,dt\\ &=r\int_0^\infty e^{-r\frac{x^p}{p}}f\left(\frac{x^p}{p}\right)x^{p-1}dx\;\; s\int_0^\infty e^{-s\frac{t^q}{q}}g\left(\frac{t^q}{q}\right)t^{q-1}\,dt.\end{aligned} \qquad (17)$$

Substituting $u=\frac{x^p}{p}, v=\frac{t^q}{q}$, $du=x^{p-1}dx$ and $dv=t^{q-1}dt$ in Equation (17) and simplifying, we obtain

$$\begin{aligned}\mathcal{G}_x^p\mathcal{G}_t^q\left[f\left(\frac{x^p}{p}\right)g\left(\frac{t^q}{q}\right)\right] &=r\int_0^\infty e^{-r\,u}\,[f(u)]du\;\; s\int_0^\infty e^{-s\,v}\,g(v)\,dv\\ &=\mathcal{G}_x[f(x)]\,\mathcal{G}_t[g(t)].\end{aligned}$$

## 4.1 CDARAT of Some Elementary Functions

In this section, we present the conformable Double ARA Transform for some basic functions.

i. Let $h\left(\frac{x^p}{p},\frac{t^q}{q}\right)=1$ . Then

$$\mathcal{G}_x^p\mathcal{G}_t^q[1]=rs\int_0^\infty\int_0^\infty e^{-r\frac{x^p}{p}-s\frac{t^q}{q}}x^{p-1}t^{q-1}dx\,dt.$$

From Property 2 and Equation (6), we get

$$\mathcal{G}_x^p\mathcal{G}_t^q[1]=\mathcal{G}_x[1]\,\mathcal{G}_t[g(1)]=1.$$

ii. Let $h\left(\frac{x^p}{p},\frac{t^q}{q}\right)=\left(\frac{x^p}{p}\right)^n\left(\frac{t^q}{q}\right)^m$. Then

$$\begin{aligned}&\mathcal{G}_x^p\mathcal{G}_t^q\left[\left(\frac{x^p}{p}\right)^n\left(\frac{t^q}{q}\right)^m\right]\\ &=rs\int_0^\infty\int_0^\infty e^{-r\frac{x^p}{p}-s\frac{t^q}{q}}\left[\left(\frac{x^p}{p}\right)^n\left(\frac{t^q}{q}\right)^m\right]x^{p-1}t^{q-1}dx\,dt\,.\end{aligned}$$

From Property 2 and Equation (6), we get

$$\mathcal{G}_x^p\mathcal{G}_t^q\left[\left(\frac{x^p}{p}\right)^n\left(\frac{t^q}{q}\right)^m\right]=\mathcal{G}_x[x^n]\,\mathcal{G}_t[t^m]=\frac{n!\,m!}{r^ns^m}.$$

iii. Let $h\left(\frac{x^p}{p},\frac{t^q}{q}\right) = e^{\lambda\frac{x^p}{p}+\beta\frac{t^q}{q}}$. Then

$$\mathcal{G}_x^p\mathcal{G}_t^q\left[e^{\lambda\frac{x^p}{p}+\beta\frac{t^q}{q}}\right] = rs\int_0^\infty\int_0^\infty e^{-r\frac{x^p}{p}-s\frac{t^q}{q}}\left[e^{\lambda\frac{x^p}{p}+\beta\frac{t^q}{q}}\right]x^{p-1}t^{q-1}dx\,dt.$$

From Property 2 and Equation (7), we get

$$\mathcal{G}_x^p\mathcal{G}_t^q\left[e^{\lambda\frac{x^p}{p}+\beta\frac{t^q}{q}}\right] = \mathcal{G}_x[e^{\lambda x}]\,\mathcal{G}_t[e^{\beta t}] = \frac{rs}{(r-\lambda)(s-\beta)}.$$

iv. Let $h\left(\frac{x^p}{p},\frac{t^q}{q}\right) = \sin\lambda\left(\frac{x^p}{p}\right)\sin\beta\left(\frac{t^q}{q}\right)$. Then

$$\mathcal{G}_x^p\mathcal{G}_t^q\left[\sin\lambda\left(\frac{x^p}{p}\right)\sin\beta\left(\frac{t^q}{q}\right)\right] = rs\int_0^\infty\int_0^\infty e^{-r\frac{x^p}{p}-s\frac{t^q}{q}}\sin\lambda\left(\frac{x^p}{p}\right)\sin\beta\left(\frac{t^q}{q}\right)x^{p-1}t^{q-1}dx\,dt.$$

From Property 2 and Equation (8), we get

$$\mathcal{G}_x^p\mathcal{G}_t^q\left[\sin\lambda\left(\frac{x^p}{p}\right)\sin\beta\left(\frac{t^q}{q}\right)\right] = \mathcal{G}_x[\sin\lambda x]\,\mathcal{G}_t[\sin\beta t] = \frac{rs}{(r^2+\lambda^2)(s^2+\lambda^2)}.$$

### 4.2 Existence Condition for the Conformable Double ARA Transform

If $h\left(\frac{x^p}{p},\frac{t^q}{q}\right)$ is an exponential order $a$ and $b$ as $\frac{x^p}{p}\to\infty, \frac{t^q}{q}\to\infty$, If $\exists$ a constant $K>0$, such that for all $x > X$ and $t > T$

$$\left|h\left(\frac{x^p}{p},\frac{t^q}{q}\right)\right| \le K\,e^{\gamma\frac{x^p}{p}+\tau\frac{t^q}{q}}.$$

it is easy to get,

$$h\left(\frac{x^p}{p},\frac{t^q}{q}\right) = O\left(e^{a\frac{x^p}{p}+b\frac{t^q}{q}}\right) \text{ as } \frac{x^p}{p}\to\infty, \frac{t^q}{q}\to\infty.$$

**Theorem 2.** Let the function $h\left(\frac{x^p}{p},\frac{t^q}{q}\right)$ be continuous on the region $(0,X)\times(0,T)$ and are of exponential orders $\gamma$ and , then the conformable double ARA transform of $h\left(\frac{x^p}{p},\frac{t^q}{q}\right)$ exists for all $Re(r) > \gamma, Re(s) > \tau$.

**Proof.** Using the definition of the CDARAT of $h\left(\frac{x^p}{p},\frac{t^q}{q}\right)$, we have

$$|G_{p.q}(r,s)| = \left|rs\int_0^\infty\int_0^\infty e^{-r\frac{x^p}{p}-s\frac{t^q}{q}}h\left(\frac{x^p}{p},\frac{t^q}{q}\right)x^{p-1}t^{q-1}dx\,dt\right|$$
$$\le rs\int_0^\infty\int_0^\infty e^{-r\frac{x^p}{p}-s\frac{t^q}{q}}\left|h\left(\frac{x^p}{p},\frac{t^q}{q}\right)\right|x^{p-1}t^{q-1}dx\,dt$$
$$\le rs\,K\int_0^\infty\int_0^\infty e^{-(r-\gamma)\frac{x^p}{p}-(s-\tau)\frac{t^q}{q}}x^{p-1}t^{q-1}dx\,dt$$
$$= \frac{rsK}{(r-\gamma)(s-\tau)}$$

For $\mathrm{Re}(r) > \gamma, \mathrm{Re}(s) > \tau$.

**Theorem 3.** Let $G_{p.q}(r,s) = \mathcal{G}_x^p\mathcal{G}_t^q\left[h\left(\frac{x^p}{p},\frac{t^q}{q}\right)\right]$, then

i. $\mathcal{G}_x^p\mathcal{G}_t^q\left[\frac{x^p}{p}h\left(\frac{x^p}{p},\frac{t^q}{q}\right)\right] = -r\frac{\partial}{\partial r}\left(\frac{1}{r}\mathcal{G}_x^p\mathcal{G}_t^q\left[h\left(\frac{x^p}{p},\frac{t^q}{q}\right)\right]\right)$.

ii. $\mathcal{G}_x^p\mathcal{G}_t^q\left[\frac{t^q}{q}h\left(\frac{x^p}{p},\frac{t^q}{q}\right)\right] = -s\frac{\partial}{\partial s}\left(\frac{1}{s}\mathcal{G}_x^p\mathcal{G}_t^q\left[h\left(\frac{x^p}{p},\frac{t^q}{q}\right)\right]\right)$.

iii. $\mathcal{G}_x^p\mathcal{G}_t^q\left[\left(\frac{x^p}{p}\right)^2h\left(\frac{x^p}{p},\frac{t^q}{q}\right)\right] = \frac{\partial^2G_{p.q}(r,s)}{\partial r^2} - \frac{2}{r}\frac{\partial G_{p.q}(r,s)}{\partial r} + \frac{2}{r^2}G_{p.q}(r,s)$.

iv. $\mathcal{G}_x^p\mathcal{G}_t^q\left[\left(\frac{t^q}{q}\right)^2h\left(\frac{x^p}{p},\frac{t^q}{q}\right)\right] = \frac{\partial^2G_{p.q}(r,s)}{\partial s^2} - \frac{2}{s}\frac{\partial G_{p.q}(r,s)}{\partial s} + \frac{2}{s^2}G_{p.q}(r,s)$.

v. $\mathcal{G}_x^p\mathcal{G}_t^q\left[\frac{x^p}{p}\frac{t^q}{q}h\left(\frac{x^p}{p},\frac{t^q}{q}\right)\right] = \frac{\partial^2G_{p.q}(r,s)}{\partial r\partial s} - \frac{1}{s}\frac{\partial G_{p.q}(r,s)}{\partial r} - \frac{1}{r}\frac{\partial G_{p.q}(r,s)}{\partial s} + \frac{1}{rs}G_{p.q}(r,s)$.

**Proof.**

**Proof of (i).** Using the definition of CDARAT of $h\left(\frac{x^p}{p},\frac{t^q}{q}\right)$, we get

$$G_{p.q}(r,s) = rs\int_0^\infty\int_0^\infty e^{-r\frac{x^p}{p}-s\frac{t^q}{q}}h\left(\frac{x^p}{p},\frac{t^q}{q}\right)x^{p-1}t^{q-1}dx\,dt. \quad (18)$$

By differentiating both sides with respect to $r$ in Equation (18), we have

$$\frac{\partial G_{p.q}(r,s)}{\partial r} = s\int_0^\infty e^{-s\frac{t^q}{q}}t^{q-1}dt\int_0^\infty\frac{\partial}{\partial r}\left(re^{-r\frac{x^p}{p}}\right)h\left(\frac{x^p}{p},\frac{t^q}{q}\right)x^{p-1}dx.$$

Calculating the partial derivative of the second integral, we can get

$$\frac{\partial G_{p.q}(r,s)}{\partial r} = s\int_0^\infty e^{-s\frac{t^q}{q}}t^{q-1}dt\int_0^\infty\left(-r\frac{x^p}{p}+1\right)e^{-r\frac{x^p}{p}}h\left(\frac{x^p}{p},\frac{t^q}{q}\right)x^{p-1}dx. \quad (19)$$

Therefore,

$$\mathcal{G}_x^p\mathcal{G}_t^q\left[\frac{x^p}{p}h\left(\frac{x^p}{p},\frac{t^q}{q}\right)\right] = -\frac{\partial G_{p.q}(r,s)}{\partial r}+\frac{1}{r}G_{p.q}(r,s) = -r\frac{\partial}{\partial r}\left(\frac{1}{r}\mathcal{G}_x^p\mathcal{G}_t^q\left[h\left(\frac{x^p}{p},\frac{t^q}{q}\right)\right]\right). \quad (20)$$

**Proof of (iii)**. Differentiating the both sides with respect to $r$ in Equation (19), we have

$$\frac{\partial^2 G_{p.q}(r,s)}{\partial r^2} = s\int_0^\infty e^{-s\frac{t^q}{q}}t^{q-1}dt\int_0^\infty \frac{\partial}{\partial r}\left[\left(1-r\frac{x^p}{p}\right)e^{-r\frac{x^p}{p}}\right]h\left(\frac{x^p}{p},\frac{t^q}{q}\right)x^{p-1}dx = s\int_0^\infty e^{-s\frac{t^q}{q}}t^{q-1}dt\int_0^\infty\left(r\left(\frac{x^p}{p}\right)^2-2\frac{x^p}{p}\right)e^{-r\frac{x^p}{p}}h\left(\frac{x^p}{p},\frac{t^q}{q}\right)x^{p-1}dx.$$

Thus,

$$\mathcal{G}_x^p\mathcal{G}_t^q\left[\left(\frac{x^p}{p}\right)^2 h\left(\frac{x^p}{p},\frac{t^q}{q}\right)\right] = \frac{\partial^2 G_{p.q}(r,s)}{\partial r^2}+\frac{2}{r}\mathcal{G}_x^p\mathcal{G}_t^q\left[\frac{x^p}{p}h\left(\frac{x^p}{p},\frac{t^q}{q}\right)\right]. \quad (21)$$

From Equation (20), we have

$$\mathcal{G}_x^p\mathcal{G}_t^q\left[\left(\frac{x^p}{p}\right)^2 h\left(\frac{x^p}{p},\frac{t^q}{q}\right)\right] = \frac{\partial^2 G_{p.q}(r,s)}{\partial r^2}+\frac{2}{r^2}G_{p.q}(r,s)-\frac{2}{r}\frac{\partial G_{p.q}(r,s)}{\partial r}.$$

Similarly, we can easily prove that:

$$\mathcal{G}_x^p\mathcal{G}_t^q\left[\frac{t^q}{q}h\left(\frac{x^p}{p},\frac{t^q}{q}\right)\right] = -s\frac{\partial}{\partial s}\left(\frac{1}{s}\mathcal{G}_x^p\mathcal{G}_t^q\left[h\left(\frac{x^p}{p},\frac{t^q}{q}\right)\right]\right).$$

$$\mathcal{G}_x^p\mathcal{G}_t^q\left[\left(\frac{t^q}{q}\right)^2 h\left(\frac{x^p}{p},\frac{t^q}{q}\right)\right] = \frac{\partial^2 G_{p.q}(r,s)}{\partial s^2}+\frac{2}{s^2}G_{p.q}(r,s)-\frac{2}{s}\frac{\partial G_{p.q}(r,s)}{\partial s}.$$

**Proof of (v)**. Differentiating both sides with respect to $s$ in Equation (19), we have

$$\frac{\partial^2 G_{p.q}(r,s)}{\partial r\partial s} = \int_0^\infty\left(-s\frac{t^q}{q}+1\right)e^{-s\frac{t^q}{q}}\left(\int_0^\infty\left(-r\frac{x^p}{p}+1\right)e^{-r\frac{x^p}{p}}h\left(\frac{x^p}{p},\frac{t^q}{q}\right)x^{p-1}dx\right)t^{q-1}dt.$$

Therefore,

$$\frac{\partial^2 G_{p.q}(r,s)}{\partial r\partial s} = \int_0^\infty s\frac{t^q}{q}e^{-s\frac{t^q}{q}}\left(\int_0^\infty r\frac{x^p}{p}e^{-r\frac{x^p}{p}}h\left(\frac{x^p}{p},\frac{t^q}{q}\right)x^{p-1}dx\right)t^{q-1}dt - \int_0^\infty e^{-s\frac{t^q}{q}}\left(\int_0^\infty r\frac{x^p}{p}e^{-r\frac{x^p}{p}}h\left(\frac{x^p}{p},\frac{t^q}{q}\right)x^{p-1}dx\right)t^{q-1}dt - \int_0^\infty s\frac{t^q}{q}e^{-s\frac{t^q}{q}}\left(\int_0^\infty e^{-r\frac{x^p}{p}}h\left(\frac{x^p}{p},\frac{t^q}{q}\right)x^{p-1}dx\right)t^{q-1}dt + \int_0^\infty e^{-s\frac{t^q}{q}}\left[\int_0^\infty e^{-r\frac{x^p}{p}}h\left(\frac{x^p}{p},\frac{t^q}{q}\right)x^{p-1}dx\right]t^{q-1}dt.$$

From (i) and (ii), we have

$$\mathcal{G}_x^p\mathcal{G}_t^q\left[\frac{x^p}{p}\frac{t^q}{q}h\left(\frac{x^p}{p},\frac{t^q}{q}\right)\right] = \frac{\partial^2 G_{p.q}(r,s)}{\partial r\partial s}+\frac{1}{rs}G_{p.q}(r,s)-\frac{1}{s}\frac{\partial G_{p.q}(r,s)}{\partial r}-\frac{1}{r}\frac{\partial G_{p.q}(r,s)}{\partial s}.$$

The proof of (ii) and (iv) can be obtained by similar arguments of (i) and (iii).

**Theorem 4.** Let $G_{p.q}(r,s)=\mathcal{G}_x^p\mathcal{G}_t^q\left[h\left(\frac{x^p}{p},\frac{t^q}{q}\right)\right]$, then

i. $\mathcal{G}_x^p\mathcal{G}_t^q\left[\frac{\partial^p}{\partial x^p}h\left(\frac{x^p}{p},\frac{t^q}{q}\right)\right] = r\,G_{p.q}(r,s) - r\,\mathcal{G}_t^q\left[h\left(0,\frac{t^q}{q}\right)\right]$.

ii. $\mathcal{G}_x^p\mathcal{G}_t^q\left[\frac{\partial^{2p}}{\partial x^{2p}}h\left(\frac{x^p}{p},\frac{t^q}{q}\right)\right] = r^2\,G_{p.q}(r,s) - r^2\,\mathcal{G}_t^q\left[h\left(0,\frac{t^q}{q}\right)\right] - r\,\mathcal{G}_t^q\left[\frac{\partial^p}{\partial x^p}h\left(0,\frac{t^q}{q}\right)\right]$.

iii. $\mathcal{G}_x^p\mathcal{G}_t^q\left[\frac{\partial^q}{\partial t^q}h\left(\frac{x^p}{p},\frac{t^q}{q}\right)\right] = s\,G_{p.q}(r,s) - s\,\mathcal{G}_x^p\left[h\left(\frac{x^p}{p},0\right)\right]$.

iv. $\mathcal{G}_x^p\mathcal{G}_t^q\left[\frac{\partial^{2q}}{\partial t^{2q}}h\left(\frac{x^p}{p},\frac{t^q}{q}\right)\right] = s^2\,G_{p.q}(r,s) - s^2\,\mathcal{G}_x^p\left[h\left(\frac{x^p}{p},0\right)\right] - s\,\mathcal{G}_x^p\left[\frac{\partial^q}{\partial t^q}h\left(\frac{x^p}{p},0\right)\right]$.

**Proof.**

i. Using the definition of CDARAT for $\frac{\partial^p}{\partial x^p}h\left(\frac{x^p}{p},\frac{t^q}{q}\right)$, we have

$$\mathcal{G}_x^p\mathcal{G}_t^q\left[\frac{\partial^p}{\partial x^p}h\left(\frac{x^p}{p},\frac{t^q}{q}\right)\right] = rs\int_0^\infty\int_0^\infty e^{-r\frac{x^p}{p}-s\frac{t^q}{q}}\frac{\partial^p}{\partial x^p}h\left(\frac{x^p}{p},\frac{t^q}{q}\right)x^{p-1}t^{q-1}dxdt \quad (22)$$

$$= s \int_0^\infty e^{-s\frac{t^q}{q}} t^{q-1} \left( r \int_0^\infty e^{-r\frac{x^p}{p}} \frac{\partial^p}{\partial x^p} h\left(\frac{x^p}{p}, \frac{t^q}{q}\right) x^{p-1} dx \right) dt.$$

Applying Property 1, $\frac{\partial^p}{\partial x^p} h\left(\frac{x^p}{p}, \frac{t^q}{q}\right) = x^{1-p} \frac{\partial}{\partial x} h\left(\frac{x^p}{p}, \frac{t^q}{q}\right)$, then Equation (22) becomes

$$\mathcal{G}_x^p \mathcal{G}_t^q \left[\frac{\partial^p}{\partial x^p} h\left(\frac{x^p}{p}, \frac{t^q}{q}\right)\right] = s \int_0^\infty e^{-s\frac{t^q}{q}} t^{q-1} \left( r \int_0^\infty e^{-r\frac{x^p}{p}} \frac{\partial}{\partial x} h(x,t) dx \right). \quad (23)$$

Thus, the integral inside bracket is given by

$$r \int_0^\infty e^{-r\frac{x^p}{p}} \frac{\partial}{\partial x} h\left(\frac{x^p}{p}, \frac{t^q}{q}\right) dx = r\left( -h\left(0, \frac{t^q}{q}\right) + r \int_0^\infty e^{-r\frac{x^p}{p}} h\left(\frac{x^p}{p}, \frac{t^q}{q}\right) x^{p-1} dx \right). \quad (24)$$

Substituting Equation (24) into Equation (23), we obtain

$$\mathcal{G}_x^p \mathcal{G}_t^q \left[\frac{\partial^p}{\partial x^p} h\left(\frac{x^p}{p}, \frac{t^q}{q}\right)\right] = r\, G_{p.q}(r,s) - r\, \mathcal{G}_t^q \left[h\left(0, \frac{t^q}{q}\right)\right]. \quad (25)$$

In the same manner, the CDARAT of $\frac{\partial^q}{\partial t^q} h\left(\frac{x^p}{p}, \frac{t^q}{q}\right)$, $\frac{\partial^{2p}}{\partial x^{2p}} h\left(\frac{x^p}{p}, \frac{t^q}{q}\right)$ and $\frac{\partial^{2q}}{\partial t^{2q}} h\left(\frac{x^p}{p}, \frac{t^q}{q}\right)$ can be obtained.

**Theorem 5.** Let $G_{p.q}(r,s) = \mathcal{G}_x^p \mathcal{G}_t^q \left[h\left(\frac{x^p}{p}, \frac{t^q}{q}\right)\right]$, then

i. $\mathcal{G}_x^p \mathcal{G}_t^q \left[\frac{x^p}{p} \frac{\partial^q}{\partial t^q} h\left(\frac{x^p}{p}, \frac{t^q}{q}\right)\right] = -rs \frac{\partial}{\partial r}\left(\frac{1}{r} \mathcal{G}_x^p \mathcal{G}_t^q \left[h\left(\frac{x^p}{p}, \frac{t^q}{q}\right)\right]\right) + rs \frac{d}{dr}\left(\frac{1}{r} \mathcal{G}_x^p \left[h\left(\frac{x^p}{p}, 0\right)\right]\right)$.

ii. $\mathcal{G}_x^p \mathcal{G}_t^q \left[\frac{t^q}{q} \frac{\partial^p}{\partial x^p} h\left(\frac{x^p}{p}, \frac{t^q}{q}\right)\right] = -rs \frac{\partial}{\partial r}\left(\frac{1}{s} \mathcal{G}_x^p \mathcal{G}_t^q \left[h\left(\frac{x^p}{p}, \frac{t^q}{q}\right)\right]\right) + rs \frac{d}{ds}\left(\frac{1}{s} \mathcal{G}_t^q \left[h\left(0, \frac{t^q}{q}\right)\right]\right)$.

**Proof of i.** The conformable double ARA transforms definition for fractional partial derivatives, implies

$$\frac{\partial}{\partial r}\left[\mathcal{G}_x^p \mathcal{G}_t^q \left(\frac{\partial^q}{\partial t^q} h\left(\frac{x^p}{p}, \frac{t^q}{q}\right)\right)\right] = s \int_0^\infty e^{-s\frac{t^q}{q}} \frac{\partial^q}{\partial t^q} h\left(\frac{x^p}{p}, \frac{t^q}{q}\right) t^{q-1} dt \left( \int_0^\infty \frac{\partial}{\partial r}\left[r\, e^{-r\frac{x^p}{p}}\right] x^{p-1} dx \right). \quad (26)$$

we calculate the partial derivative in the second integral as follows

$$\int_0^\infty \frac{\partial}{\partial r}\left(r\, e^{-r\frac{x^p}{p}}\right) x^{p-1} dx = \int_0^\infty e^{-r\frac{x^p}{p}} x^{p-1} dx - r \int_0^\infty \frac{x^p}{p} e^{-r\frac{x^p}{p}} x^{p-1} dx. \quad (27)$$

Substituting Equation (27) into Equation (26), we get

$$\frac{\partial}{\partial r}\left[\mathcal{G}_x^p \mathcal{G}_t^q \left(\frac{\partial^q}{\partial t^q} h\left(\frac{x^p}{p}, \frac{t^q}{q}\right)\right)\right] = s \int_0^\infty e^{-s\frac{t^q}{q}} \int_0^\infty e^{-r\frac{x^p}{p}} \frac{\partial^q}{\partial t^q} h\left(\frac{x^p}{p}, \frac{t^q}{q}\right) x^{p-1} t^{q-1} dx dt - rs \int_0^\infty e^{-s\frac{t^q}{q}} \int_0^\infty e^{-r\frac{x^p}{p}} \frac{x^p}{p} \frac{\partial^q h\left(\frac{x^p}{p}, \frac{t^q}{q}\right)}{\partial t^q} x^{p-1} t^{q-1} dx\, dt.$$

Thus,

$$\mathcal{G}_x^p \mathcal{G}_t^q \left[\frac{x^p}{p} \frac{\partial^q}{\partial t^q} h\left(\frac{x^p}{p}, \frac{t^q}{q}\right)\right] = -\frac{\partial}{\partial r} \mathcal{G}_x^p \mathcal{G}_t^q \left[\frac{\partial^q}{\partial t^q} h\left(\frac{x^p}{p}, \frac{t^q}{q}\right)\right] + \frac{1}{r} \mathcal{G}_x^p \mathcal{G}_t^q \left[\frac{\partial^q}{\partial t^q} h\left(\frac{x^p}{p}, \frac{t^q}{q}\right)\right].$$

Using Theorem 4, we have

$$\mathcal{G}_x^p \mathcal{G}_t^q \left[\frac{x^p}{p} \frac{\partial^q h\left(\frac{x^p}{p}, \frac{t^q}{q}\right)}{\partial t^q}\right] = -rs \frac{\partial}{\partial r}\left(\frac{1}{r} \mathcal{G}_x^p \mathcal{G}_t^q \left[h\left(\frac{x^p}{p}, \frac{t^q}{q}\right)\right]\right) + rs \frac{d}{dr}\left(\frac{1}{r} \mathcal{G}_x^p \left[h\left(\frac{x^p}{p}, 0\right)\right]\right)$$

Similarly, one can prove that

$$\mathcal{G}_x^p \mathcal{G}_t^q \left[\frac{t^q}{q} \frac{\partial^p}{\partial x^p} h\left(\frac{x^p}{p}, \frac{t^q}{q}\right)\right] = -rs \frac{\partial}{\partial s}\left(\frac{1}{s} \mathcal{G}_x^p \mathcal{G}_t^q \left[h\left(\frac{x^p}{p}, \frac{t^q}{q}\right)\right]\right) + rs \frac{d}{ds}\left(\frac{1}{s} \mathcal{G}_t^q \left[h\left(0, \frac{t^q}{q}\right)\right]\right).$$

In the following, we introduce the previous results in the following table, Table 1, below:

Table 1. Analysis of the presented results

| $h\left(\frac{x^p}{p},\frac{t^q}{q}\right)$ | $\mathcal{G}_x^p\mathcal{G}_t^q\left[h\left(\frac{x^p}{p},\frac{t^q}{q}\right)\right]=G_{p.q}(r,s)$ |
|---|---|
| 1 | 1 |
| $\left(\frac{x^p}{p}\right)^n\left(\frac{t^q}{q}\right)^m$ | $\frac{n!\,m!}{r^n s^m}.$ |
| $e^{\lambda\frac{x^p}{p}+\beta\frac{t^q}{q}}$ | $\frac{rs}{(r-\lambda)(s-\beta)}.$ |
| $\sin\lambda\left(\frac{x^p}{p}\right)\sin\beta\left(\frac{t^q}{q}\right)$ | $\frac{rs}{(r^2+\lambda^2)(s^2+\lambda^2)}.$ |
| $\frac{x^p}{p}h\left(\frac{x^p}{p},\frac{t^q}{q}\right)$ | $-r\frac{\partial}{\partial r}\left(\frac{1}{r}\mathcal{G}_x^p\mathcal{G}_t^q\left[h\left(\frac{x^p}{p},\frac{t^q}{q}\right)\right]\right).$ |
| $\frac{t^q}{q}h\left(\frac{x^p}{p},\frac{t^q}{q}\right)$ | $-s\frac{\partial}{\partial s}\left(\frac{1}{s}\mathcal{G}_x^p\mathcal{G}_t^q\left[h\left(\frac{x^p}{p},\frac{t^q}{q}\right)\right]\right).$ |
| $\left(\frac{x^p}{p}\right)^2h\left(\frac{x^p}{p},\frac{t^q}{q}\right)$ | $\frac{\partial^2G_{p.q}(r,s)}{\partial r^2}-\frac{2}{r}\frac{\partial G_{p.q}(r,s)}{\partial r}+\frac{2}{r^2}G_{p.q}(r,s).$ |
| $\left(\frac{t^q}{q}\right)^2h\left(\frac{x^p}{p},\frac{t^q}{q}\right)$ | $\frac{\partial^2G_{p.q}(r,s)}{\partial s^2}-\frac{2}{s}\frac{\partial G_{p.q}(r,s)}{\partial s}+\frac{2}{s^2}G_{p.q}(r,s).$ |
| $\frac{x^p}{p}\frac{t^q}{q}h\left(\frac{x^p}{p},\frac{t^q}{q}\right)$ | $\frac{\partial^2G_{p.q}(r,s)}{\partial r\partial s}-\frac{1}{s}\frac{\partial G_{p.q}(r,s)}{\partial r}-\frac{1}{r}\frac{\partial G_{p.q}(r,s)}{\partial s}+\frac{1}{rs}G_{p.q}(r,s).$ |
| $\frac{\partial^p}{\partial x^p}h\left(\frac{x^p}{p},\frac{t^q}{q}\right)$ | $r\,G_{p.q}(r,s)-r\,\mathcal{G}_t^q\left[h\left(0,\frac{t^q}{q}\right)\right].$ |
| $\frac{\partial^q}{\partial t^q}h\left(\frac{x^p}{p},\frac{t^q}{q}\right)$ | $s\,G_{p.q}(r,s)-s\,\mathcal{G}_x^p\left[h\left(\frac{x^p}{p},0\right)\right].$ |
| $\frac{\partial^{2p}}{\partial x^{2p}}h\left(\frac{x^p}{p},\frac{t^q}{q}\right)$ | $r^2\,G_{p.q}(r,s)-r^2\,\mathcal{G}_t^q\left[h\left(0,\frac{t^q}{q}\right)\right]-r\,\mathcal{G}_t^q\left[\frac{\partial^p}{\partial x^p}h\left(0,\frac{t^q}{q}\right)\right]..$ |
| $\frac{\partial^{2q}}{\partial t^{2q}}h\left(\frac{x^p}{p},\frac{t^q}{q}\right)$ | $s^2\,G_{p.q}(r,s)-s^2\,\mathcal{G}_x^p\left[h\left(\frac{x^p}{p},0\right)\right]-s\,\mathcal{G}_x^p\left[\frac{\partial^q}{\partial t^q}h\left(\frac{x^p}{p},0\right)\right].$ |
| $\frac{x^p}{p}\frac{\partial^q}{\partial t^q}h\left(\frac{x^p}{p},\frac{t^q}{q}\right)$ | $-rs\frac{\partial}{\partial r}\left(\frac{1}{r}\mathcal{G}_x^p\mathcal{G}_t^q\left[h\left(\frac{x^p}{p},\frac{t^q}{q}\right)\right]\right)+rs\frac{d}{dr}\left(\frac{1}{r}\mathcal{G}_x^p\left[h\left(\frac{x^p}{p},0\right)\right]\right).$ |
| $\frac{t^q}{q}\frac{\partial^p}{\partial x^p}h\left(\frac{x^p}{p},\frac{t^q}{q}\right)$ | $-rs\frac{\partial}{\partial s}\left(\frac{1}{s}\mathcal{G}_x^p\mathcal{G}_t^q\left[h\left(\frac{x^p}{p},\frac{t^q}{q}\right)\right]\right)+rs\frac{d}{ds}\left(\frac{1}{s}\mathcal{G}_t^q\left[h\left(0,\frac{t^q}{q}\right)\right]\right).$ |

## 5 Applications

The CDARADM is used in this section of the study to solve regular and singular one-dimensional conformable fractional coupled Burger's equations. The goal problem is the same as the problem examined in [1], when $p=1$ and $q=1$. This is what we mention here.

**Example 1.**

Consider the One-dimensional conformable fractional coupled Burgers' equation of the form

$$\begin{aligned}\frac{\partial^q u}{\partial t^q}-\frac{\partial^{2p}u}{\partial x^{2p}}+\lambda\,u\frac{\partial^p u}{\partial x^p}+\alpha\,\frac{\partial^p}{\partial x^p}(uv)&=k\left(\frac{x^p}{p},\frac{t^q}{q}\right)\\ \frac{\partial^q v}{\partial t^q}-\frac{\partial^{2p}v}{\partial x^{2p}}+\lambda\,v\frac{\partial^p v}{\partial x^p}+\beta\,\frac{\partial^p}{\partial x^p}(uv)&=l\left(\frac{x^p}{p},\frac{t^q}{q}\right),\end{aligned}\tag{28}$$

subject to

$$\begin{aligned}u\left(\frac{x^p}{p},0\right)&=k_1\left(\frac{x^p}{p}\right),\\ v\left(\frac{x^p}{p},0\right)&=l_1\left(\frac{x^p}{p}\right),\end{aligned}\tag{29}$$

for $t>0$. Here, $k\left(\frac{x^p}{p},\frac{t^q}{q}\right)$, $l\left(\frac{x^p}{p},\frac{t^q}{q}\right)$, $k_1\left(\frac{x^p}{p}\right)$ and $l_1\left(\frac{x^p}{p}\right)$ are given functions, $\lambda,\alpha,$ and $\beta$ are arbitrary parameters depending on the Peclet number, Stokes velocity of particles due to gravity and Brownian diffusivity, see, [9]. Now, operating the conformable double ARA transform to Equation (28) and the single conformable single ARA transform for Equation (29), to get

$$U(r,s)=K_1(r)+\frac{K(r,s)}{s}+\frac{1}{s}\,\mathcal{G}_x^p\mathcal{G}_t^q\left[\frac{\partial^{2p}u}{\partial x^{2p}}-\lambda\,u\frac{\partial^p u}{\partial x^p}-\alpha\,\frac{\partial^p}{\partial x^p}(uv)\right].\tag{30}$$

$$V(r,s)=L_1(r)+\frac{L(r,s)}{s}+\frac{1}{s}\,\mathcal{G}_x^p\mathcal{G}_t^q\left[\frac{\partial^{2p}v}{\partial x^{2p}}-\lambda\,v\frac{\partial^p v}{\partial x^p}-\beta\,\frac{\partial^p}{\partial x^p}(uv)\right].\tag{31}$$

The CDARADM defines the solution of the target problem $u\left(\frac{x^p}{p},\frac{t^q}{q}\right)$ and $v\left(\frac{x^p}{p},\frac{t^q}{q}\right)$ in the form of infinite series as

$$\begin{aligned}u\left(\frac{x^p}{p},\frac{t^q}{q}\right)&=\sum_{n=0}^{\infty}u_n\left(\frac{x^p}{p},\frac{t^q}{q}\right),\\ v\left(\frac{x^p}{p},\frac{t^q}{q}\right)&=\sum_{n=0}^{\infty}v_n\left(\frac{x^p}{p},\frac{t^q}{q}\right).\end{aligned}\tag{32}$$

Define the Adomian's polynomials $A_n$, $B_n$ and $C_n$ as

$$A_n=\sum_{n=0}^{\infty}u_n u_{xn},\tag{33}$$

$$B_n = \sum_{n=0}^{\infty} v_n v_{xn},$$
$$C_n = \sum_{n=0}^{\infty} u_n v_{n.}$$

We can compute the Adomian polynomials of the nonlinear terms $uu_x, vv_x$ and $uv$ by the formulas

$$A_0 = u_0 u_{0x.}$$
$$A_1 = u_0 u_{1x} + u_1 u_{0x.}$$
$$A_2 = u_0 u_{2x} + u_1 u_{1x} + u_2 u_{0x}.$$
$$A_3 = u_0 u_{3x} + u_1 u_{2x} + u_2 u_{1x} + u_3 u_{0x}.$$
$$A_3 = u_0 u_{4x} + u_1 u_{3x} + u_2 u_{2x} + u_3 u_{1x} + u_4 u_{0x.}$$
$$\vdots$$
$$B_0 = v_0 v_{0x.}$$
$$B_1 = v_0 v_{1x} + v_1 v_{0x}.$$
$$B_2 = v_0 v_{2x} + v_1 v_{1x} + v_2 v_{0x}.$$
$$B_3 = v_0 v_{3x} + v_1 v_{2x} + v_2 v_{1x} + v_3 v_{0x}.$$
$$B_3 = v_0 v_{4x} + v_1 v_{3x} + v_2 v_{2x} + v_3 v_{1x} + v_4 v_{0x.}$$
$$\vdots$$
$$C_0 = u_0 v_{0.}$$
$$C_1 = u_0 v_1 + u_1 v_0.$$
$$C_2 = u_0 v_2 + u_1 v_1 + u_2 v_0.$$
$$C_3 = u_0 v_3 + u_1 v_2 + u_2 v_1 + u_3 v_0.$$
$$C_3 = u_0 v_4 + u_1 v_3 + u_2 v_2 + u_3 v_1 + u_4 v_{0.}$$
$$\vdots$$

Operating the inverse double ARA transform to Equation (30) and Equation (31), utilizing Equation (33), we get

$$\sum_{n=0}^{\infty} u_n\left(\frac{x^p}{p},\frac{t^q}{q}\right) = k_1(x) \; \mathcal{G}_r^{-1}\mathcal{G}_s^{-1}\left[\frac{K(r,s)}{s}\right] + \mathcal{G}_r^{-1}\mathcal{G}_s^{-1}\left[\frac{1}{s}\,\mathcal{G}_x^p\mathcal{G}_t^q\left(\frac{\partial^{2p}u_n}{\partial x^{2p}}\right)\right] - \mathcal{G}_r^{-1}\mathcal{G}_s^{-1}\left[\frac{1}{s}\,\mathcal{G}_x^p\mathcal{G}_t^q(\lambda\, A_n)\right] - \mathcal{G}_r^{-1}\mathcal{G}_s^{-1}\left[\frac{1}{s}\,\mathcal{G}_x^p\mathcal{G}_t^q(\alpha\, C_n)\right], \tag{34}$$

and

$$\sum_{n=0}^{\infty} v_n\left(\frac{x^p}{p},\frac{t^q}{q}\right) = l_1(x) + \mathcal{G}_r^{-1}\mathcal{G}_s^{-1}\left[\frac{L(r,s)}{s}\right] + \mathcal{G}_r^{-1}\mathcal{G}_s^{-1}\left[\frac{1}{s}\,\mathcal{G}_x^p\mathcal{G}_t^q\left(\frac{\partial^{2p}v_n}{\partial x^{2p}}\right)\right] - \mathcal{G}_r^{-1}\mathcal{G}_s^{-1}\left[\frac{1}{s}\,\mathcal{G}_x^p\mathcal{G}_t^q(\lambda\, B_n)\right] - \mathcal{G}_r^{-1}\mathcal{G}_s^{-1}\left[\frac{1}{s}\,\mathcal{G}_x^p\mathcal{G}_t^q(\beta\, C_n)\right]. \tag{35}$$

Now, we compare both sides of Equation (34) and Equation (35), to get

$$u_0 = k_1(x) + \mathcal{G}_r^{-1}\mathcal{G}_s^{-1}\left[\frac{F(r,s)}{s}\right],$$
$$v_0 = l_1(x) + \mathcal{G}_r^{-1}\mathcal{G}_s^{-1}\left[\frac{L(r,s)}{s}\right]. \tag{36}$$

Following that, the recursive relation can be expressed as

$$u_n = \mathcal{G}_r^{-1}\mathcal{G}_s^{-1}\left[\frac{1}{s}\,\mathcal{G}_x^p\mathcal{G}_t^q\left(\frac{\partial^{2p}v_n}{\partial x^{2\alpha}}\right)\right] - \mathcal{G}_r^{-1}\mathcal{G}_s^{-1}\left[\frac{1}{s}\,\mathcal{G}_x^p\mathcal{G}_t^q(\lambda\, A_n)\right] - \mathcal{G}_r^{-1}\mathcal{G}_s^{-1}\left[\frac{1}{s}\,\mathcal{G}_x^p\mathcal{G}_t^q(\alpha\, C_n)\right], \tag{37}$$

and

$$v_n = \mathcal{G}_r^{-1}\mathcal{G}_s^{-1}\left[\frac{1}{s}\,\mathcal{G}_x^p\mathcal{G}_t^q\left(\frac{\partial^{2p}u_n}{\partial x^{2\alpha}}\right)\right] - \mathcal{G}_r^{-1}\mathcal{G}_s^{-1}\left[\frac{1}{s}\,\mathcal{G}_x^p\mathcal{G}_t^q(\lambda\, B_n)\right] - \mathcal{G}_r^{-1}\mathcal{G}_s^{-1}\left[\frac{1}{s}\,\mathcal{G}_x^p\mathcal{G}_t^q(\beta\, C_n)\right]. \tag{38}$$

Herein, we should state that the solutions in (37) and (38) exist, provided the inverse double ARA transform exists $\forall p$ and $s$.

Putting $\lambda = -2$, $\alpha = \beta = 1$ and $k\left(\frac{x^p}{p},\frac{t^q}{q}\right) = l\left(\frac{x^p}{p},\frac{t^q}{q}\right) = 0$ in Equation (28) and $k_1\left(\frac{x^p}{p}\right) = l_1\left(\frac{x^p}{p}\right) = \sin\left(\frac{x^p}{p}\right)$ in Equation (29), we obtain the one-dimensional homogeneous coupled Burgers fractional equation is the conformable sense

$$\frac{\partial^q u}{\partial t^q} - \frac{\partial^{2p}u}{\partial x^{2p}} - 2u\frac{\partial^p u}{\partial x^p} + \frac{\partial^p}{\partial x^p}(uv) = 0,$$
$$\frac{\partial^q v}{\partial t^q} - \frac{\partial^{2p}v}{\partial x^{2p}} - 2v\frac{\partial^p v}{\partial x^p} + \frac{\partial^p}{\partial x^p}(uv) = 0. \tag{39}$$

with initial condition

$$u\left(\frac{x^p}{p},0\right) = \sin\left(\frac{x^p}{p}\right),$$
$$v\left(\frac{x^p}{p},0\right) = \sin\left(\frac{x^p}{p}\right). \tag{40}$$

By using Equations (36)–(38) , we have

$$u_0 = \sin\left(\frac{x^p}{p}\right), \qquad v_0 = \sin\left(\frac{x^p}{p}\right),$$

$$u_1 = \mathcal{G}_r^{-1}\mathcal{G}_s^{-1}\left[\frac{1}{s}\,\mathcal{G}_x^p\mathcal{G}_t^q\left[\frac{\partial^{2p}u_0}{\partial x^{2p}} + 2u_0\frac{\partial^p u_0}{\partial x^p} - \frac{\partial^p}{\partial x^p}(u_0 v_0)\right]\right]$$
$$= \mathcal{G}_r^{-1}\mathcal{G}_s^{-1}\left[\frac{1}{s}\,\mathcal{G}_x^p\mathcal{G}_t^q\left[-\sin\left(\frac{x^p}{p}\right)\right]\right]$$
$$= \mathcal{G}_r^{-1}\mathcal{G}_s^{-1}\left[\frac{-r}{s(r^2+1)}\right]$$
$$= -\frac{t^q}{q}\sin\left(\frac{x^p}{p}\right),$$

$$v_1 = \mathcal{G}_r^{-1}\mathcal{G}_s^{-1}\left[\frac{1}{s}\,\mathcal{G}_x^p\mathcal{G}_t^q\left[\frac{\partial^{2p}v_0}{\partial x^{2p}} + 2v_0\frac{\partial^p v_0}{\partial x^p} - \frac{\partial^p}{\partial x^p}(u_0 v_0)\right]\right]$$
$$= \mathcal{G}_r^{-1}\mathcal{G}_s^{-1}\left[\frac{1}{s}\,\mathcal{G}_x^p\mathcal{G}_t^q\left[-\sin\left(\frac{x^p}{p}\right)\right]\right]$$
$$= \mathcal{G}_r^{-1}\mathcal{G}_s^{-1}\left[\frac{-r}{s(r^2+1)}\right]$$
$$= -\frac{t^q}{q}\sin\left(\frac{x^p}{p}\right),$$

$$u_2 = \mathcal{G}_r^{-1}\mathcal{G}_s^{-1}\left[\frac{1}{s}\,\mathcal{G}_x^p\mathcal{G}_t^q\left[\frac{\partial^{2p}u_1}{\partial x^{2p}} + 2\left(u_0\frac{\partial^p u_1}{\partial x^p} + u_1\frac{\partial^p u_0}{\partial x^p}\right) - \frac{\partial^p}{\partial x^p}(u_0 v_1 + u_1 v_0)\right]\right]$$
$$= \mathcal{G}_r^{-1}\mathcal{G}_s^{-1}\left[\frac{1}{s}\,\mathcal{G}_x^p\mathcal{G}_t^q\left[\frac{t^q}{q}\sin\left(\frac{x^p}{p}\right)\right]\right]$$
$$= \mathcal{G}_r^{-1}\mathcal{G}_s^{-1}\left[\frac{r}{s^2(r^2+1)}\right]$$
$$= \frac{\left(\frac{t^q}{q}\right)^2}{2}\sin\left(\frac{x^p}{p}\right),$$

$$v_2 = \mathcal{G}_r^{-1}\mathcal{G}_s^{-1}\left[\frac{1}{s}\mathcal{G}_x^p\mathcal{G}_t^q\left[\frac{\partial^{2p}v_1}{\partial x^{2p}} + 2\left(v_0\frac{\partial^p v_1}{\partial x^p} + v_1\frac{\partial^p v_0}{\partial x^p}\right) - \frac{\partial^p}{\partial x^p}(u_0 v_1 + u_1 v_0)\right]\right]$$
$$= \mathcal{G}_r^{-1}\mathcal{G}_s^{-1}\left[\frac{1}{s}\mathcal{G}_x^p\mathcal{G}_t^q\left[\frac{t^q}{q}\sin\left(\frac{x^p}{p}\right)\right]\right]$$
$$= \mathcal{G}_r^{-1}\mathcal{G}_s^{-1}\left[\frac{r}{s^2(r^2+1)}\right]$$
$$= \frac{\left(\frac{t^q}{q}\right)^2}{2}\sin\left(\frac{x^p}{p}\right),$$

and

$$u_3 = \mathcal{G}_r^{-1}\mathcal{G}_s^{-1}\left[\frac{1}{s}\,\mathcal{G}_x^p\mathcal{G}_t^q\left[\frac{\partial^{2p}u_2}{\partial x^{2p}} + 2\left(u_0\frac{\partial^p u_2}{\partial x^p} + u_1\frac{\partial^p u_1}{\partial x^p} + u_2\frac{\partial^p u_0}{\partial x^p}\right) - \frac{\partial^p}{\partial x^p}(u_0 v_2 + u_1 v_1 + u_2 v_0)\right]\right]$$
$$= -\frac{\left(\frac{t^q}{q}\right)^3}{6}\sin\left(\frac{x^p}{p}\right),$$

$$v_3 = \mathcal{G}_r^{-1}\mathcal{G}_s^{-1}\left[\frac{1}{s}\,\mathcal{G}_x^p\mathcal{G}_t^q\left[\frac{\partial^{2p}v_2}{\partial x^{2p}} + 2\left(v_0\frac{\partial^p v_2}{\partial x^p} + v_1\frac{\partial^p v_1}{\partial x^p} + v_2\frac{\partial^p v_0}{\partial x^p}\right) - \frac{\partial^p}{\partial x^p}(u_0 v_2 + u_1 v_1 + u_2 v_0)\right]\right]$$
$$= -\frac{\left(\frac{t^q}{q}\right)^3}{6}\,sin\left(\frac{x^p}{p}\right).$$

Therefore, using Equation (32), we can express the series solution as

$$u\left(\frac{x^p}{p},\frac{t^q}{q}\right)=u_0+u_1+u_2+\cdots$$
$$=\left(1-\left(\frac{t^q}{q}\right)+\frac{\left(\frac{t^q}{q}\right)^2}{2!}-\frac{\left(\frac{t^q}{q}\right)^3}{3!}+\cdots\right)\sin\left(\frac{x^p}{p}\right),$$
$$\left(\frac{x^p}{p},\frac{t^q}{q}\right)=v_0+v_1+v_2+\cdots$$
$$=\left(1-\left(\frac{t^q}{q}\right)+\frac{\left(\frac{t^q}{q}\right)^2}{2!}-\frac{\left(\frac{t^q}{q}\right)^3}{6!}+\cdots\right)v\sin\left(\frac{x^p}{p}\right),$$

and hence the exact solutions become

$$u\left(\frac{x^p}{p},\frac{t^q}{q}\right)=e^{-\frac{t^q}{q}}\sin\left(\frac{x^p}{p}\right),$$
$$v\left(\frac{x^p}{p},\frac{t^q}{q}\right)=e^{-\frac{t^q}{q}}\sin\left(\frac{x^p}{p}\right).$$

By taking $p=1$ and $q=1$, the fractional solution of Equation (39) becomes

$$u\left(\frac{x^p}{p},\frac{t^q}{q}\right)=e^{-t}\sin x,$$
$$v\left(\frac{x^p}{p},\frac{t^q}{q}\right)=e^{-t}\sin x.$$

The behavior of the velocity field of the two-CDARADM (28) and (29) is depicted in Figure 1 for (a) the approximate and exact solutions of $u\left(\frac{x^p}{p},\frac{t^q}{q}\right)$ for Example 1, when $p=q$, at $p=0.8,0.9,1$, and (b) the approximate and exact solutions of $u\left(\frac{x^p}{p},\frac{t^q}{q}\right)$, for Example 1, when taking various values of fractional order $q$ ($q=0.8,0.9,1$) and $p=1$.

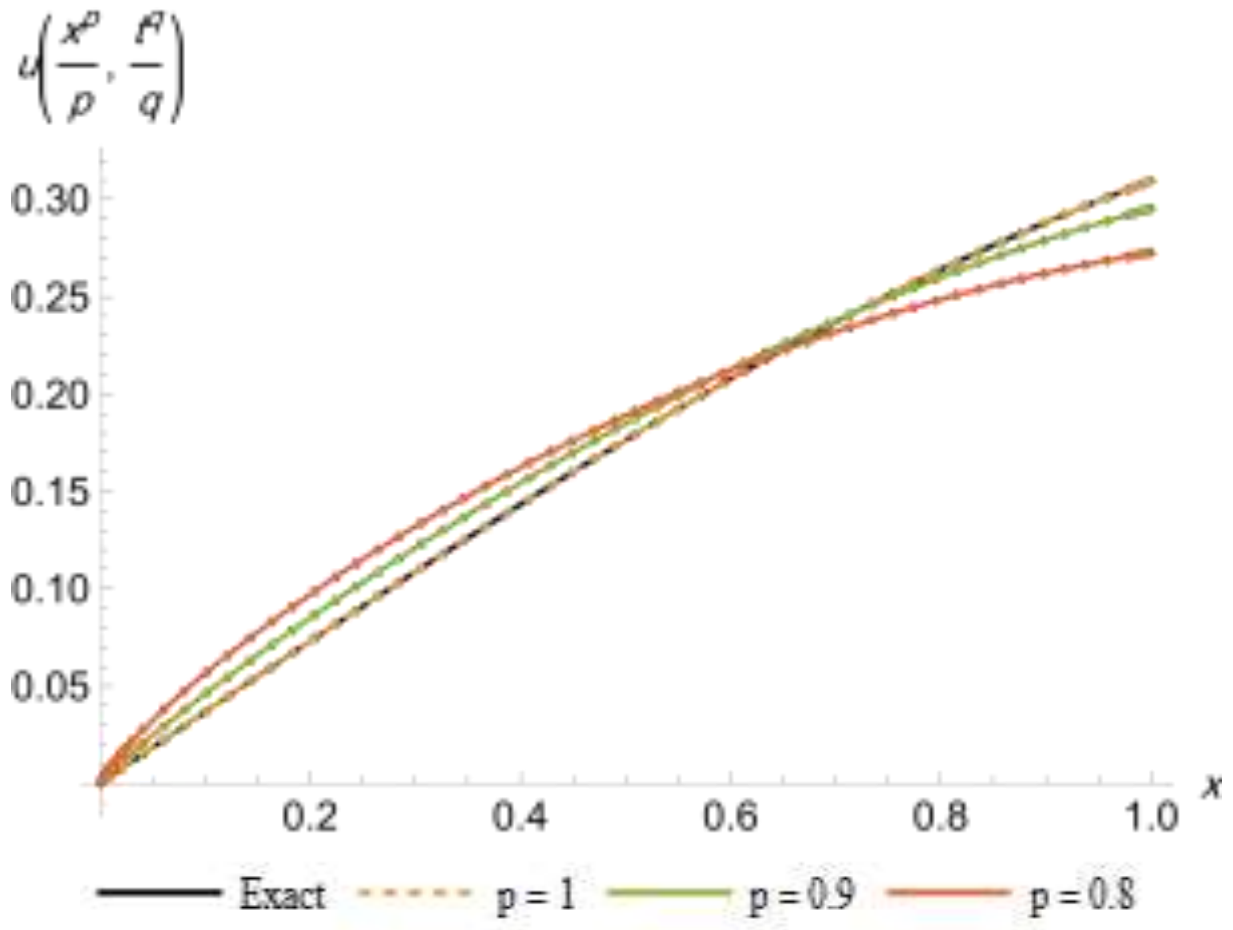

(a)

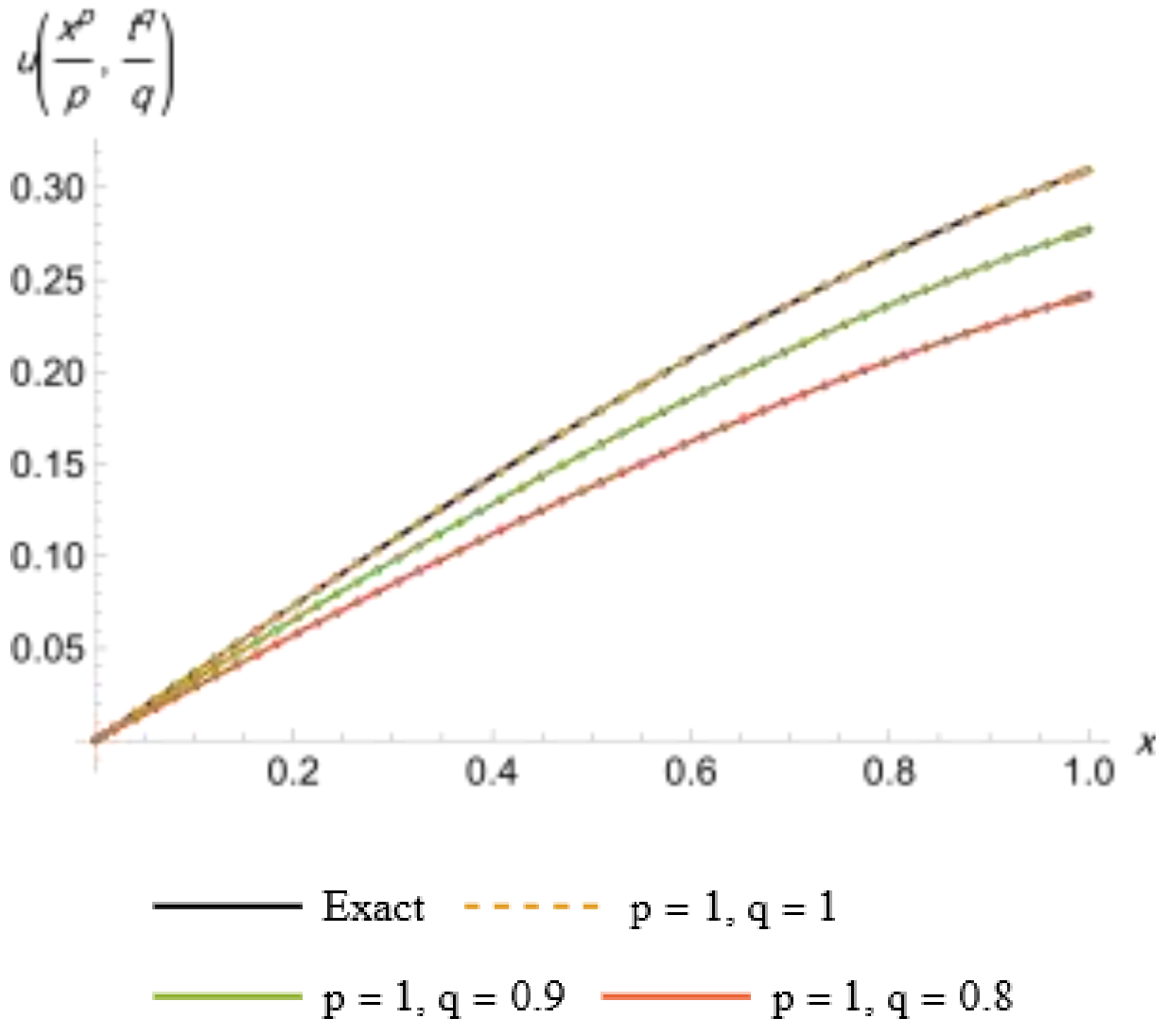

(b)

Fig. 1: The figure of (a) The approximate and exact solutions of $u(x^p/p, t^q/q)$ for Example 5.1, when $p=q$, (b) The approximate and exact solutions of $u(x^p/p, t^q/q)$, for Example 1, when taking various values of $q$ and $p=1$.

Table 2 below, presents the absolute errors considering $p=q=1$, $x=1, t\in[0.1,0.5]$. We compare the obtained results by exact solution.

Table 2. Analysis of error to $u\left(\frac{x^p}{p},\frac{t^q}{q}\right)$ for Example 1, on $x=1, t\in[0.1,0.5]$

| | Exact Solution | Approximate Solution | $\left\|u\left(\frac{x^p}{p},\frac{t^q}{q}\right)-u_6\left(\frac{x^p}{p},\frac{t^q}{q}\right)\right\|$ |
|---|---|---|---|
| 0.1 | 0.761394 | 0.761394 | $1.64894\times10^{-11}$ |
| 0.15 | 0.724261 | 0.724261 | $2.80003\times10^{-10}$ |
| 0.2 | 0.688938 | 0.688938 | $2.08481\times10^{-9}$ |
| 0.25 | 0.655338 | 0.655338 | $9.88052\times10^{-9}$ |
| 0.3 | 0.623377 | 0.623377 | $3.51889\times10^{-8}$ |
| 0.35 | 0.592975 | 0.592975 | $1.02897\times10^{-7}$ |
| 0.4 | 0.564055 | 0.564055 | $2.60452\times10^{-7}$ |
| 0.45 | 0.536546 | 0.536546 | $5.90459\times10^{-7}$ |
| 0.5 | 0.510378 | 0.510379 | $1.22715\times10^{-6}$ |

**Example 2.**
Consider the singular one-dimensional conformable fractional coupled Burgers' equation with the Bessel operator of the form

$$\begin{aligned}&\frac{\partial^q u}{\partial t^q}-\frac{p}{x^p}\frac{\partial^p}{\partial x^p}\left(\frac{x^p}{p}\frac{\partial^p}{\partial x^p}u\right)+\lambda\,u\frac{\partial^p u}{\partial x^p}+\alpha\,\frac{\partial^p}{\partial x^p}(uv)=k\left(\frac{x^p}{p},\frac{t^q}{q}\right),\\&\frac{\partial^q v}{\partial t^q}-\frac{p}{x^p}\frac{\partial^p}{\partial x^p}\left(\frac{x^p}{p}\frac{\partial^p}{\partial x^p}v\right)+\lambda\,v\frac{\partial^p v}{\partial x^p}+\beta\,\frac{\partial^p}{\partial x^p}(uv)=l\left(\frac{x^p}{p},\frac{t^q}{q}\right),\end{aligned}\tag{41}$$

and with initial conditions

$$\begin{aligned}u\left(\frac{x^p}{p},0\right)&=\ k_1\left(\frac{x^p}{p}\right),\\v\left(\frac{x^p}{p},0\right)&=\ l_1\left(\frac{x^p}{p}\right).\end{aligned}\tag{42}$$

Where the linear terms $\frac{p}{x^p}\frac{\partial^p}{\partial x^p}\left(\frac{x^p}{p}\frac{\partial^p}{\partial x^p}\right)$ is known as a conformable Bessel operator where $\lambda,\alpha$, and $\beta$ are real constants.

Multiplying both sides of Equation (41) by $\frac{x^p}{p}$, we get

$$\begin{aligned}&\frac{x^p}{p}\frac{\partial^q u}{\partial t^q}-\frac{\partial^p}{\partial x^p}\left(\frac{x^p}{p}\frac{\partial^p}{\partial x^p}u\right)+\lambda\,\frac{x^p}{p}u\frac{\partial^p u}{\partial x^p}+\alpha\,\frac{x^p}{p}\frac{\partial^p}{\partial x^p}(uv)=\frac{x^p}{p}k\left(\frac{x^p}{p},\frac{t^q}{q}\right),\\&\frac{x^p}{p}\frac{\partial^q v}{\partial t^q}-\frac{\partial^p}{\partial x^p}\left(\frac{x^p}{p}\frac{\partial^p}{\partial x^p}v\right)+\lambda\,\frac{x^p}{p}v\frac{\partial^p v}{\partial x^p}+\beta\,\frac{x^p}{p}\frac{\partial^p}{\partial x^p}(uv)=\frac{x^p}{p}l\left(\frac{x^p}{p},\frac{t^q}{q}\right).\end{aligned}\tag{43}$$

Applying conformable double ARA transform to both sides of Equation (43) and single conformable ARA transform for the initial condition, we get

$$\begin{aligned}&\mathcal{G}_x^p\mathcal{G}_t^q\left[\frac{x^p}{p}\frac{\partial^q u}{\partial t^q}\right]-\mathcal{G}_x^p\mathcal{G}_t^q\left[\frac{\partial^p}{\partial x^p}\left(\frac{x^p}{p}\frac{\partial^p}{\partial x^p}u\right)-\lambda\,\frac{x^p}{p}u\frac{\partial^p u}{\partial x^p}-\alpha\,\frac{x^p}{p}\frac{\partial^p}{\partial x^p}(uv)\right]=\mathcal{G}_x^p\mathcal{G}_t^q\left[\frac{x^p}{p}k\left(\frac{x^p}{p},\frac{t^q}{q}\right)\right],\\&\mathcal{G}_x^p\mathcal{G}_t^q\left[\frac{x^p}{p}\frac{\partial^q v}{\partial t^q}\right]-\mathcal{G}_x^p\mathcal{G}_t^q\left[\frac{\partial^p}{\partial x^p}\left(\frac{x^p}{p}\frac{\partial^p}{\partial x^p}v\right)-\lambda\,\frac{x^p}{p}u\frac{\partial^p v}{\partial x^p}-\beta\,\frac{x^p}{p}\frac{\partial^p}{\partial x^p}(uv)\right]=\mathcal{G}_x^p\mathcal{G}_t^q\left[\frac{x^p}{p}l\left(\frac{x^p}{p},\frac{t^q}{q}\right)\right].\end{aligned}\tag{44}$$

Applying Theorem 3 and Theorem 5, we have

$$\begin{aligned}&-rs\frac{\partial}{\partial r}\left(\frac{1}{r}U(r,s)\right)+rs\frac{d}{dr}\left(\frac{1}{r}\mathcal{G}_x^p[k_1(x)]\right)=\mathcal{G}_x^p\mathcal{G}_t^q\left[\frac{\partial^p}{\partial x^p}\left(\frac{x^p}{p}\frac{\partial^p}{\partial x^p}u\right)-\lambda\,\frac{x^p}{p}u\frac{\partial^p u}{\partial x^p}-\alpha\,\frac{x^p}{p}\frac{\partial^p}{\partial x^p}(uv)\right]-r\frac{\partial}{\partial r}\left(\frac{1}{r}\mathcal{G}_x^p\mathcal{G}_t^q\left[k\left(\frac{x^p}{p},\frac{t^q}{q}\right)\right]\right),\\&-rs\frac{\partial}{\partial r}\left(\frac{1}{r}V(r,s)\right)+rs\frac{d}{dr}\left(\frac{1}{r}\mathcal{G}_x^p[l_1(x)]\right)=\mathcal{G}_x^p\mathcal{G}_t^q\left[\frac{\partial^p}{\partial x^p}\left(\frac{x^p}{p}\frac{\partial^p}{\partial x^p}v\right)-\lambda\,\frac{x^p}{p}v\frac{\partial^p u}{\partial x^p}-\beta\,\frac{x^p}{p}\frac{\partial^p}{\partial x^p}(uv)\right]-r\frac{\partial}{\partial r}\left(\frac{1}{r}\mathcal{G}_x^p\mathcal{G}_t^q\left[l\left(\frac{x^p}{p},\frac{t^q}{q}\right)\right]\right).\end{aligned}\tag{45}$$

Simplifying Equation (45), we obtain

$$\begin{aligned}
&\frac{\partial}{\partial r}\left(\frac{1}{r}U(r,s)\right)\\
&=\frac{d}{dr}\left(\frac{1}{r}\mathcal{G}_x^p[k_1(x)]\right)\\
&-\frac{1}{rs}\mathcal{G}_x^p\mathcal{G}_t^q\left[\frac{\partial^p}{\partial x^p}\left(\frac{x^p}{p}\frac{\partial^p}{\partial x^p}u\right)\right.\\
&\left.-\ \lambda\ \frac{x^p}{p}u\frac{\partial^p u}{\partial x^p}-\alpha\ \frac{x^p}{p}\frac{\partial^p}{\partial x^p}(uv)\right]\\
&+\frac{1}{s}\frac{\partial}{\partial r}\left(\frac{1}{r}\mathcal{G}_x^p\mathcal{G}_t^q\left[k\left(\frac{x^p}{p},\frac{t^q}{q}\right)\right]\right),\\
&\frac{\partial}{\partial r}\left(\frac{1}{r}V(r,s)\right)\\
&=\frac{d}{dr}\left(\frac{1}{r}\mathcal{G}_x^p[l_1(x)]\right)\\
&-\frac{1}{rs}\mathcal{G}_x^p\mathcal{G}_t^q\left[\frac{\partial^p}{\partial x^p}\left(\frac{x^p}{p}\frac{\partial^p}{\partial x^p}v\right)\right.\\
&\left.-\ \lambda\ \frac{x^p}{p}v\frac{\partial^p u}{\partial x^p}-\beta\frac{x^p}{p}\frac{\partial^p}{\partial x^p}(uv)\right]\\
&+\frac{1}{s}\frac{\partial}{\partial r}\left(\frac{1}{r}\mathcal{G}_x^p\mathcal{G}_t^q\left[l\left(\frac{x^p}{p},\frac{t^q}{q}\right)\right]\right).
\end{aligned}\tag{46}$$

Applying the definite integral $\int_0^r$ with respect to $r$ to both sides of Equation (46)

$$\begin{aligned}
&\frac{1}{r}U(r,s)\\
&=\int_0^r\frac{d}{dr}\left(\frac{1}{r}\mathcal{G}_x^p[k_1(x)]\right)dr\\
&-\frac{1}{s}\int_0^r\left(\frac{1}{r}\mathcal{G}_x^p\mathcal{G}_t^q\left[\frac{\partial^p}{\partial x^p}\left(\frac{x^p}{p}\frac{\partial^p}{\partial x^p}u\right)-\ \lambda\ \frac{x^p}{p}N_1\right.\right.\\
&\left.\left.-\ \alpha\ \frac{x^p}{p}N_2\right]\right)dr\\
&+\frac{1}{s}\int_0^r\left(\frac{\partial}{\partial r}\left(\frac{1}{r}\mathcal{G}_x^p\mathcal{G}_t^q\left[k\left(\frac{x^p}{p},\frac{t^q}{q}\right)\right]\right)\right)dr,\\
&\frac{1}{r}V(r,s)\\
&=\int_0^r\frac{d}{dr}\left(\frac{1}{r}\mathcal{G}_x^p[l_1(x)]\right)dr\\
&-\frac{1}{s}\int_0^r\left(\frac{1}{r}\mathcal{G}_x^p\mathcal{G}_t^q\left[\frac{\partial^p}{\partial x^p}\left(\frac{x^p}{p}\frac{\partial^p}{\partial x^p}v\right)-\ \lambda\ \frac{x^p}{p}N_3\right.\right.\\
&\left.\left.-\ \beta\ \frac{x^p}{p}N_2\right]\right)dr\\
&+\frac{1}{s}\int_0^r\left(\frac{\partial}{\partial r}\left(\frac{1}{r}\mathcal{G}_x^p\mathcal{G}_t^q\left[l\left(\frac{x^p}{p},\frac{t^q}{q}\right)\right]\right)\right)dr.
\end{aligned}$$

Multiplying both sides of the equations by $r$ , we get

$$\begin{aligned}
&U(r,s)\\
&=r\int_0^r\frac{d}{dr}\left(\frac{1}{r}\mathcal{G}_x^p[k_1(x)]\right)dr\\
&-\frac{r}{s}\int_0^r\left(\frac{1}{r}\mathcal{G}_x^p\mathcal{G}_t^q\left[\frac{\partial^p}{\partial x^p}\left(\frac{x^p}{p}\frac{\partial^p}{\partial x^p}u\right)\right.\right.\\
&\left.\left.-\ \lambda\ \frac{x^p}{p}N_1-\alpha\ \frac{x^p}{p}N_2\right]\right)dr\\
&+\frac{r}{s}\int_0^r\left(\frac{\partial}{\partial r}\left(\frac{1}{r}\mathcal{G}_x^p\mathcal{G}_t^q\left[k\left(\frac{x^p}{p},\frac{t^q}{q}\right)\right]\right)\right)dr,\\
&V(r,s)\\
&=r\int_0^r\frac{d}{dr}\left(\frac{1}{r}\mathcal{G}_x^p[l_1(x)]\right)dr\\
&-\frac{r}{s}\int_0^r\left(\frac{1}{r}\mathcal{G}_x^p\mathcal{G}_t^q\left[\frac{\partial^p}{\partial x^p}\left(\frac{x^p}{p}\frac{\partial^p}{\partial x^p}v\right)\right.\right.\\
&\left.\left.-\ \lambda\ \frac{x^p}{p}N_3-\beta\ \frac{x^p}{p}N_2\right]\right)dr\\
&+\frac{r}{s}\int_0^r\left(\frac{\partial}{\partial r}\left(\frac{1}{r}\mathcal{G}_x^p\mathcal{G}_t^q\left[l\left(\frac{x^p}{p},\frac{t^q}{q}\right)\right]\right)\right)dr.
\end{aligned}\tag{47}$$

Utilizing the CDARADM to present the solution of $u\left(\frac{x^p}{p},\frac{t^q}{q}\right)$ and $v\left(\frac{x^p}{p},\frac{t^q}{q}\right)$ by infinite series as

$$\begin{aligned}
u\left(\frac{x^p}{p},\frac{t^q}{q}\right)&=\sum_{n=0}^{\infty}u_n\left(\frac{x^p}{p},\frac{t^q}{q}\right),\\
v\left(\frac{x^p}{p},\frac{t^q}{q}\right)&=\sum_{n=0}^{\infty}v_n\left(\frac{x^p}{p},\frac{t^q}{q}\right).
\end{aligned}\tag{48}$$

Define the nonlinear operators as

$$\begin{aligned}
N_1=\sum_{n=0}^{\infty}A_n,\qquad N_2=\sum_{n=0}^{\infty}C_n,\\
N_3=\sum_{n=0}^{\infty}B_n.
\end{aligned}\tag{49}$$

Operating the double inverse transform to Equation (47) and making use of Equation (48) and Equation (49), we have

$$\begin{aligned}
&\sum_{n=0}^{\infty}u_n\left(\frac{x^p}{p},\frac{t^q}{q}\right)\\
&=k_1(x)+\mathcal{G}_r^{-1}\mathcal{G}_s^{-1}\left[\frac{1}{s}\left(\mathcal{G}_x^p\mathcal{G}_t^q\left[k\left(\frac{x^p}{p},\frac{t^q}{q}\right)\right]\right)\right]\\
&-\ \mathcal{G}_r^{-1}\mathcal{G}_s^{-1}\left[\frac{r}{s}\int_0^r\frac{1}{r}\left(\mathcal{G}_x^p\mathcal{G}_t^q\left[\frac{\partial^p}{\partial x^p}\left(\frac{x^p}{p}\left(\sum_{n=0}^{\infty}u_n\right)\right)\right]\right)dr\right]\\
&+\mathcal{G}_r^{-1}\mathcal{G}_s^{-1}\left[\frac{r}{s}\int_0^r\frac{1}{r}\left(\mathcal{G}_x^p\mathcal{G}_t^q\left[\lambda\ \frac{x^p}{p}\sum_{n=0}^{\infty}A_n\right]\right)dr\right]\\
&+\mathcal{G}_r^{-1}\mathcal{G}_s^{-1}\left[\frac{r}{s}\int_0^r\frac{1}{r}\left(\mathcal{G}_x^p\mathcal{G}_t^q\left[\alpha\ \frac{x^p}{p}\sum_{n=0}^{\infty}C_n\right]\right)dr\right].
\end{aligned}\tag{50}$$

And

$$\begin{aligned}&\sum_{n=0}^{\infty} v_n\left(\frac{x^p}{p},\frac{t^q}{q}\right)\\&= l_1(x) + \mathcal{G}_r^{-1}\mathcal{G}_s^{-1}\left[\frac{1}{s}\left(\mathcal{G}_x^p\mathcal{G}_t^q\left[l\left(\frac{x^p}{p},\frac{t^q}{q}\right)\right]\right)\right]\\&- \mathcal{G}_r^{-1}\mathcal{G}_s^{-1}\left[\frac{r}{s}\int_0^r\frac{1}{r}\left(\mathcal{G}_x^p\mathcal{G}_t^q\left[\frac{\partial^p}{\partial x^p}\left(\frac{x^p}{p}\left(\sum_{n=0}^{\infty} v_n\right)\right)\right]\right)dr\right]\\&+ \mathcal{G}_r^{-1}\mathcal{G}_s^{-1}\left[\frac{r}{s}\int_0^r\frac{1}{r}\left(\mathcal{G}_x^p\mathcal{G}_t^q\left[\lambda\,\frac{x^p}{p}\sum_{n=0}^{\infty} B_n\right]\right)dr\right]\\&+ \mathcal{G}_r^{-1}\mathcal{G}_s^{-1}\left[\frac{r}{s}\int_0^r\frac{1}{r}\left(\mathcal{G}_x^p\mathcal{G}_t^q\left[\beta\,\frac{x^p}{p}\sum_{n=0}^{\infty} C_n\right]\right)dr\right].\end{aligned} \tag{51}$$

Now, we can express the first few components as

$$\begin{aligned}&u_0\\&= k_1(x)\\&+ \mathcal{G}_r^{-1}\mathcal{G}_s^{-1}\left[\frac{1}{s}\left(\mathcal{G}_x^p\mathcal{G}_t^q\left[k\left(\frac{x^p}{p},\frac{t^q}{q}\right)\right]\right)\right],\\&v_0 = l_1(x)\\&+ \mathcal{G}_r^{-1}\mathcal{G}_s^{-1}\left[\frac{1}{s}\left(\mathcal{G}_x^p\mathcal{G}_t^q\left[l\left(\frac{x^p}{p},\frac{t^q}{q}\right)\right]\right)\right].\end{aligned} \tag{52}$$

And

$$\begin{aligned}&u_{n+1}\left(\frac{x^p}{p},\frac{t^q}{q}\right)\\&= - \mathcal{G}_r^{-1}\mathcal{G}_s^{-1}\left[\frac{r}{s}\int_0^r\frac{1}{r}\left(\mathcal{G}_x^p\mathcal{G}_t^q\left[\frac{\partial^r}{\partial x^r}\left(\frac{x^p}{p}\left(\sum_{n=0}^{\infty} u_n\right)\right)\right]\right)dr\right]\\&+ \mathcal{G}_r^{-1}\mathcal{G}_s^{-1}\left[\frac{r}{s}\int_0^r\frac{1}{r}\left(\mathcal{G}_x^p\mathcal{G}_t^q\left[\lambda\,\frac{x^p}{p}\sum_{n=0}^{\infty} A_n\right]\right)dr\right]\\&+ \mathcal{G}_r^{-1}\mathcal{G}_s^{-1}\left[\frac{r}{s}\int_0^r\frac{1}{r}\left(\mathcal{G}_x^p\mathcal{G}_t^q\left[\alpha\,\frac{x^p}{p}\sum_{n=0}^{\infty} C_n\right]\right)dr\right].\end{aligned} \tag{53}$$

And

$$\begin{aligned}&v_{n+1}\left(\frac{x^p}{p},\frac{t^q}{q}\right)\\&= - \mathcal{G}_r^{-1}\mathcal{G}_s^{-1}\left[\frac{r}{s}\int_0^r\frac{1}{r}\left(\mathcal{G}_x^p\mathcal{G}_t^q\left[\frac{\partial^r}{\partial x^r}\left(\frac{x^p}{p}\left(\sum_{n=0}^{\infty} v_n\right)\right)\right]\right)dr\right]\\&+ \mathcal{G}_r^{-1}\mathcal{G}_s^{-1}\left[\frac{r}{s}\int_0^r\frac{1}{r}\left(\mathcal{G}_x^p\mathcal{G}_t^q\left[\lambda\,\frac{x^p}{p}\sum_{n=0}^{\infty} B_n\right]\right)dr\right]\\&+ \mathcal{G}_r^{-1}\mathcal{G}_s^{-1}\left[r\int_0^r\frac{1}{r}\left(\mathcal{G}_x^p\mathcal{G}_t^q\left[\beta\,\frac{x^p}{p}\sum_{n=0}^{\infty} C_n\right]\right)dr\right].\end{aligned} \tag{54}$$

In addition, we assume that the double inverse transform in (53) and (54) exist, and substitution $\lambda = -2,\ \alpha = \beta = 1$ and $k\left(\frac{x^p}{p},\frac{t^q}{q}\right) = l\left(\frac{x^p}{p},\frac{t^q}{q}\right) = \left(\frac{x^p}{p}\right)^2 e^{\frac{t^q}{q}} - 4e^{\frac{t^q}{q}}$ in Equation (41) and $k_1\left(\frac{x^p}{p}\right) = l_1\left(\frac{x^p}{p}\right) = \left(\frac{x^p}{p}\right)^2$ in Equation (42), we obtain the singular conformable coupled Burgers fractional equation of one dimensional

$$\begin{aligned}&\frac{\partial^q u}{\partial t^q} - \frac{p}{x^p}\frac{\partial^p}{\partial x^p}\left(\frac{x^p}{p}\frac{\partial^p}{\partial x^p}u\right) - 2u\frac{\partial^p u}{\partial x^p}\\&\quad + \frac{\partial^p}{\partial x^p}(uv)\\&\quad = \left(\frac{x^p}{p}\right)^2 e^{\frac{t^q}{q}} - 4e^{\frac{t^q}{q}},\\&\frac{\partial^q v}{\partial t^q} - \frac{p}{x^p}\frac{\partial^p}{\partial x^p}\left(\frac{x^p}{p}\frac{\partial^p}{\partial x^p}v\right) - 2v\frac{\partial^p v}{\partial x^p}\\&\quad + \frac{\partial^p}{\partial x^p}(uv)\\&\quad = \left(\frac{x^p}{p}\right)^2 e^{\frac{t^q}{q}} - 4e^{\frac{t^q}{q}},\end{aligned} \tag{55}$$

subject to

$$\begin{aligned}u\left(\frac{x^p}{p},0\right) &= \left(\frac{x^p}{p}\right)^2,\\ v\left(\frac{x^p}{p},0\right) &= \left(\frac{x^p}{p}\right)^2.\end{aligned} \tag{56}$$

By following similar steps, we obtain

$$\begin{aligned}&\sum_{n=0}^{\infty} u_n\left(\frac{x^p}{p},\frac{t^q}{q}\right)\\&= \mathcal{G}_r^{-1}\mathcal{G}_s^{-1}\left[\frac{1}{s}\left(\mathcal{G}_x^p\mathcal{G}_t^q\left[\left(\frac{x^p}{p}\right)^2 e^{\frac{t^q}{q}} - 4e^{\frac{t^q}{q}}\right]\right)\right]\\&- \mathcal{G}_r^{-1}\mathcal{G}_s^{-1}\left[\frac{r}{s}\int_0^r\frac{1}{r}\left(\mathcal{G}_x^p\mathcal{G}_t^q\left[\frac{\partial^p}{\partial x^p}\left(\frac{x^p}{p}\left(\sum_{n=0}^{\infty} u_n\right)\right)\right]\right)dr\right]\\&- \mathcal{G}_r^{-1}\mathcal{G}_s^{-1}\left[\frac{r}{s}\int_0^r\frac{1}{r}\left(\mathcal{G}_x^p\mathcal{G}_t^q\left[2\frac{x^p}{p}\sum_{n=0}^{\infty} A_n\right]\right)dr\right]\\&+ \mathcal{G}_r^{-1}\mathcal{G}_s^{-1}\left[\frac{r}{s}\int_0^r\frac{1}{r}\left(\mathcal{G}_x^p\mathcal{G}_t^q\left[\frac{x^p}{p}\sum_{n=0}^{\infty} C_n\right]\right)dr\right] + \left(\frac{x^p}{p}\right)^2,\end{aligned} \tag{57}$$

$$\begin{aligned}&\sum_{n=0}^{\infty} v_n\left(\frac{x^p}{p},\frac{t^q}{q}\right)\\&= \mathcal{G}_r^{-1}\mathcal{G}_s^{-1}\left[\frac{1}{s}\left(\mathcal{G}_x^p\mathcal{G}_t^q\left[\left(\frac{x^p}{p}\right)^2 e^{\frac{t^q}{q}} - 4e^{\frac{t^q}{q}}\right]\right)\right]\\&- \mathcal{G}_r^{-1}\mathcal{G}_s^{-1}\left[\frac{r}{s}\int_0^r\frac{1}{r}\left(\mathcal{G}_x^p\mathcal{G}_t^q\left[\frac{\partial^p}{\partial x^p}\left(\frac{x^p}{p}\left(\sum_{n=0}^{\infty} v_n\right)\right)\right]\right)dr\right]\\&- \mathcal{G}_r^{-1}\mathcal{G}_s^{-1}\left[\frac{r}{s}\int_0^r\frac{1}{r}\left(\mathcal{G}_x^p\mathcal{G}_t^q\left[2\frac{x^p}{p}\sum_{n=0}^{\infty} B_n\right]\right)dr\right]\\&+ \mathcal{G}_r^{-1}\mathcal{G}_s^{-1}\left[\frac{r}{s}\int_0^r\frac{1}{r}\left(\mathcal{G}_x^p\mathcal{G}_t^q\left[\frac{x^p}{p}\sum_{n=0}^{\infty} C_n\right]\right)dr\right] + \left(\frac{x^p}{p}\right)^2.\end{aligned} \tag{58}$$

Using Equations (52) – (54) the components are given by

$$\begin{aligned}u_0 &= \left(\frac{x^p}{p}\right)^2 + \mathcal{G}_r^{-1}\mathcal{G}_s^{-1}\left[\frac{1}{s}\mathcal{G}_x^p\mathcal{G}_t^q\left(\left(\frac{x^p}{p}\right)^2 e^{\frac{t^q}{q}} - 4e^{\frac{t^q}{q}}\right)\right]\\&= \left(\frac{x^p}{p}\right)^2 + \mathcal{G}_r^{-1}\mathcal{G}_s^{-1}\left[\frac{1}{s}\left(\frac{2}{r^2}\left(\frac{s}{s-1}\right) - \frac{4s}{s-1}\right)\right]\end{aligned}$$

$$= \left(\frac{x^p}{p}\right)^2 + \mathcal{G}_r^{-1}\mathcal{G}_s^{-1}\left[\frac{1}{s}\left(\frac{2}{r^2}\left(\frac{s^2}{s-1} - s\right) + 4s\left(1 - \frac{s}{s-1}\right)\right)\right]$$
$$= \left(\frac{x^p}{p}\right)^2 + \mathcal{G}_r^{-1}\mathcal{G}_s^{-1}\left[\frac{1}{s}\left(\frac{2}{r^2}\left(\frac{s^2}{s-1}\right) - \frac{4s^2}{s-1} - \frac{2s}{r^2} + 4s\right)\right]$$
$$= \left(\frac{x^p}{p}\right)^2 + \mathcal{G}_r^{-1}\mathcal{G}_s^{-1}\left[\left(\frac{2}{r^2}\left(\frac{s}{s-1}\right) - \frac{4s}{s-1} - \frac{2}{r^2} + 4\right)\right]$$
$$= \left(\frac{x^p}{p}\right)^2 + \left(\frac{x^p}{p}\right)^2 e^{\frac{t^q}{q}} - 4e^{\frac{t^q}{q}} - \left(\frac{x^p}{p}\right)^2 + 4$$
$$= \left(\frac{x^p}{p}\right)^2 e^{\frac{t^q}{q}} - 4e^{\frac{t^q}{q}} + 4,$$

$u_1$
$$= -\mathcal{G}_r^{-1}\mathcal{G}_s^{-1}\left[\frac{r}{s}\int_0^r \frac{1}{r}\left(\mathcal{G}_x^p\mathcal{G}_t^q\left[\frac{\partial^p}{\partial x^p}\left(\frac{x^p}{p}\frac{\partial^p u_0}{\partial x^p}\right) + 2\frac{x^p}{p}u_0\frac{\partial^p u_0}{\partial x^p}\right.\right.\right.$$
$$\left.\left.\left. - \frac{x^p}{p}\frac{\partial^p}{\partial x^p}(u_0\, v_0)\right]\right) dr\right]$$
$$= -\mathcal{G}_r^{-1}\mathcal{G}_s^{-1}\left[\frac{r}{s}\int_0^r \frac{1}{r}\left(\mathcal{G}_x^p\mathcal{G}_t^q\left[4\frac{x^p}{p}e^{\frac{t^q}{q}}\right.\right.\right.$$
$$+ 4\left(\left(\frac{x^p}{p}\right)^2 e^{\frac{t^q}{q}} - 4e^{\frac{t^q}{q}} + 4\right)$$
$$\left.\left.\left. - 4\left(\left(\frac{x^p}{p}\right)^2 e^{\frac{t^q}{q}} - 4e^{\frac{t^q}{q}} + 4\right)\right]\right) dr\right]$$
$$= -\mathcal{G}_r^{-1}\mathcal{G}_s^{-1}\left[\frac{r}{s}\int_0^r \frac{1}{r}\left(\mathcal{G}_x^p\mathcal{G}_t^q\left[\left(4\frac{x^p}{p}e^{\frac{t^q}{q}}\right)\right]\right) dr\right]$$
$$= \mathcal{G}_r^{-1}\mathcal{G}_s^{-1}\left[\frac{1}{s}\int_0^r \left(\frac{4s}{s-1}\right) dr\right] = \mathcal{G}_r^{-1}\mathcal{G}_s^{-1}\left[\frac{4}{s}\left(\frac{s^2}{s-1} - s\right)\right]$$
$$= 4e^{\frac{t^q}{q}} - 4,$$

$u_2 = 0 = u_3 = \cdots$

In a similar way, we obtain

$$v_0 = \left(\frac{x^p}{p}\right)^2 e^{\frac{t^q}{q}} - 4e^{\frac{t^q}{q}} + 4\,.$$
$$v_1 = 4e^{\frac{t^q}{q}} - 4.$$
$$v_2 = 0 = v_3 = \cdots.$$

Following that, one can express the solution as

$$u\left(\frac{x^p}{p}, \frac{t^q}{q}\right) = u_0 + u_1,$$
$$v\left(\frac{x^p}{p}, \frac{t^q}{q}\right) = v_0 + v_1\,\cdot$$

Therefore, the exact solution is given by

$$u\left(\frac{x^p}{p}, \frac{t^q}{q}\right) = \left(\frac{x^p}{p}\right)^2 e^{\frac{t^q}{q}},$$
$$v\left(\frac{x^p}{p}, \frac{t^q}{q}\right) = \left(\frac{x^p}{p}\right)^2 e^{\frac{t^q}{q}}.$$

By taking $p = 1$ and $q = 1$, the fractional solution becomes

$$u(x,t) = x^2 e^t \quad , \quad v(x,t) = x^2 e^t.$$

The behavior of the velocity field of the two-CDARADM (41) and (42) is depicted in Figure 2 for (a) the approximate and exact solutions of $u$ for Example 2, at $p = q = 0.8, 0.9, 1$, and (b) the approximate and exact solutions of $u\left(\frac{x^p}{p}, \frac{t^q}{q}\right)$, for Example 2, when taking various values of $q$ ($q = 0.8, 0.9, 1$) and $p = 1$.

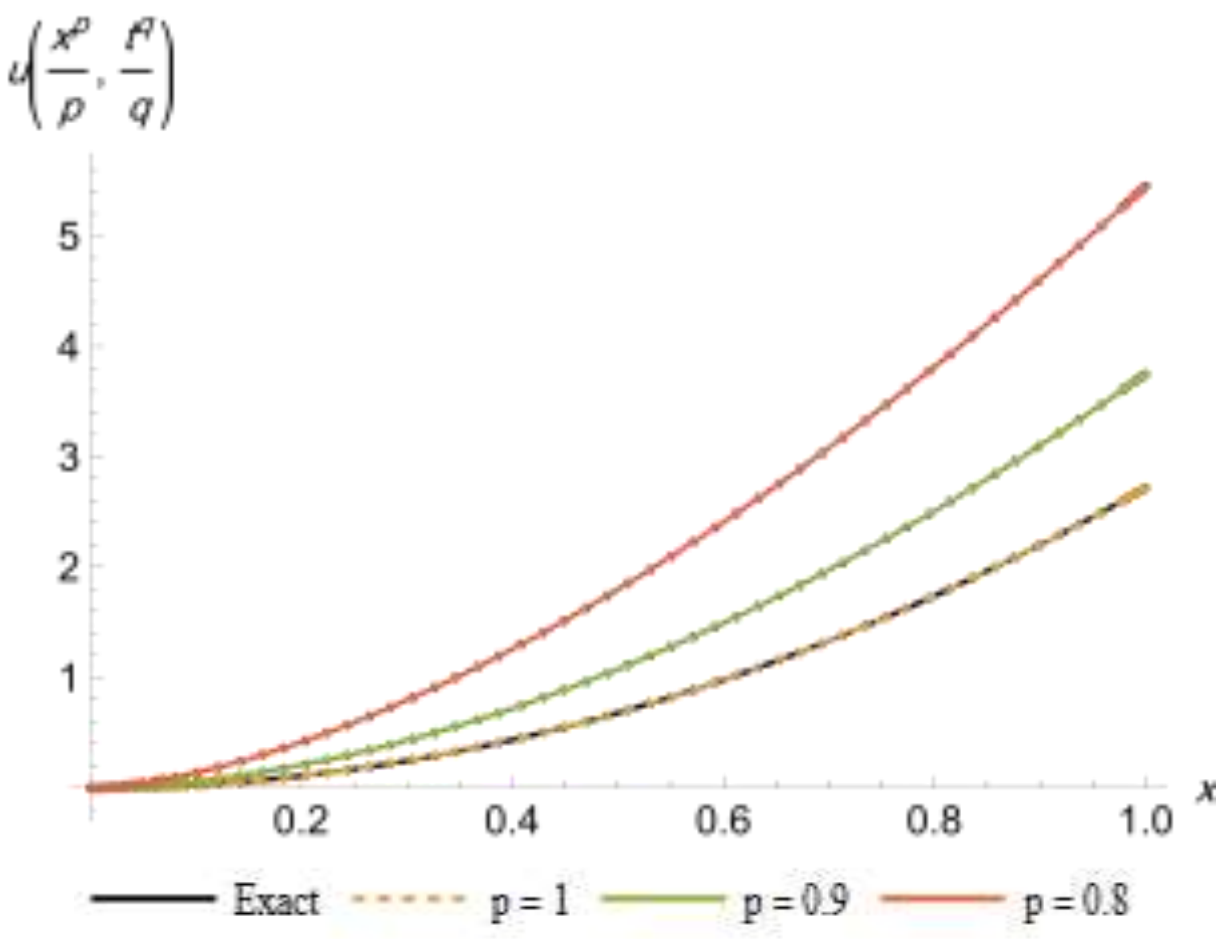

(a)

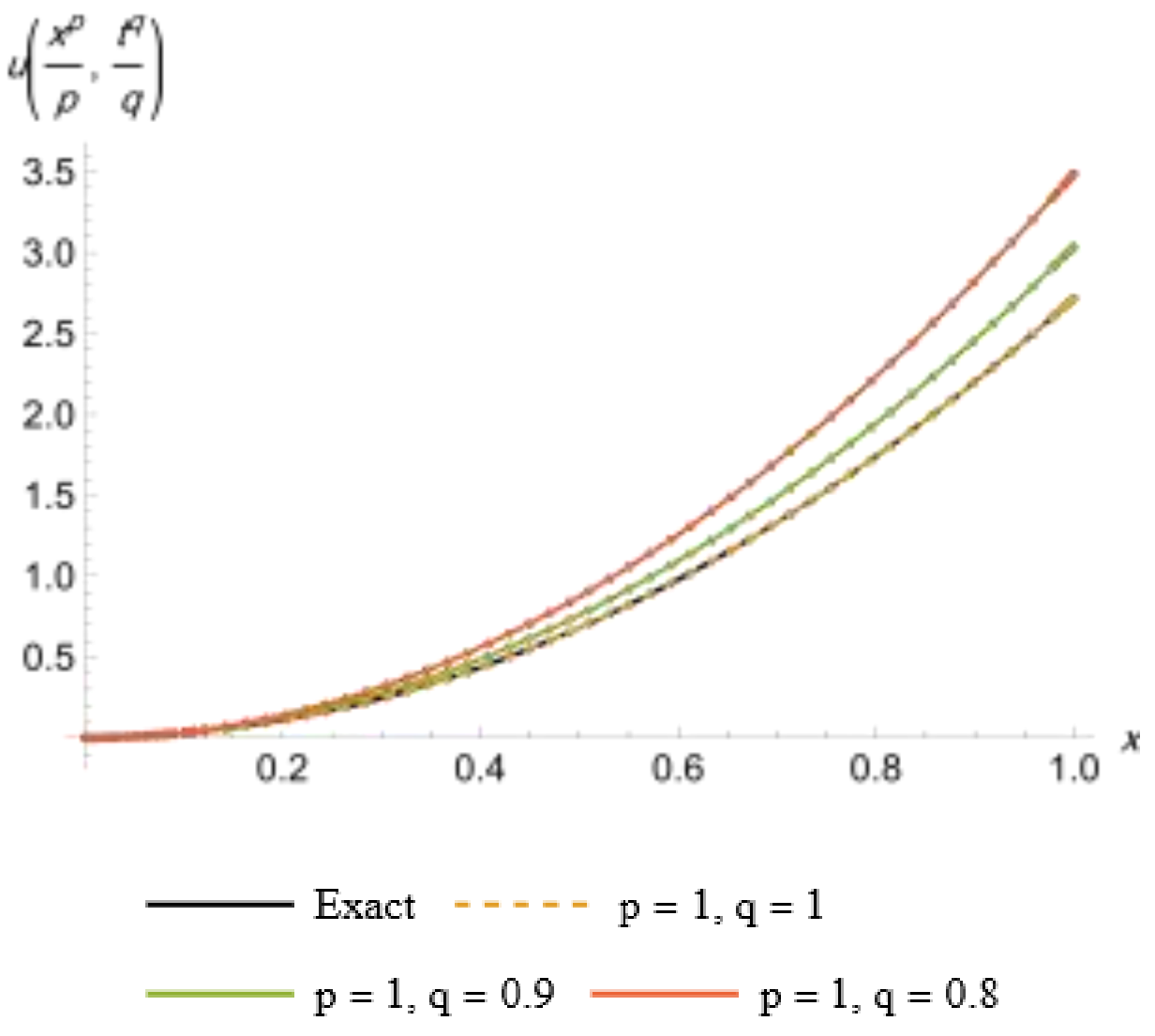

(b)

Fig. 2: The figure of (a) The approximate and exact solutions of $u\left(\frac{x^p}{p}, \frac{t^q}{q}\right)$ for Example 5.2, when $p = q$, (b) The approximate and exact solutions of $u\left(\frac{x^p}{p}, \frac{t^q}{q}\right)$, for Example 2, when we take different values of fractional order $q$ and $p = 1$.

Table 3 below, presents the absolute errors with respect to $p = q = 1$, $x = 1, t \in [0.1, 0.5]$. We compare the obtained results by exact solution.

Table 3. Error analysis of $u\left(\frac{x^p}{p},\frac{t^q}{q}\right)$ for Example 2 on $x=1, t\in[0.1,0.5]$

| | **Exact Solution** | **Approximate Solution** | $\left\|u\left(\frac{x^p}{p},\frac{t^q}{q}\right)-u_6\left(\frac{x^p}{p},\frac{t^q}{q}\right)\right\|$ |
|---|---|---|---|
| 0.1 | 1.10517 | 1.10517 | $2.00924\times10^{-11}$ |
| 0.15 | 1.16183 | 1.16183 | $3.45471\times10^{-10}$ |
| 0.2 | 1.2214 | 1.2214 | $2.60461\times10^{-9}$ |
| 0.25 | 1.28403 | 1.28403 | $1.24994\times10^{-8}$ |
| 0.3 | 1.34986 | 1.34986 | $4.5076\times10^{-8}$ |
| 0.35 | 1.41907 | 1.41907 | $1.33467\times10^{-7}$ |
| 0.4 | 1.49182 | 1.49182 | $3.42086\times10^{-7}$ |
| 0.45 | 1.56831 | 1.56831 | $7.85295\times10^{-7}$ |
| 0.5 | 1.64872 | 1.64872 | $1.65264\times10^{-6}$ |

# 6 Conclusion

In the current study, we defined and went over some of the characteristics of the conformable double ARA transform. The conformable double ARA decomposition method is a novel approach that we present for the solution of nonlinear conformable partial differential equations. We used the proposed approach, a novel amalgamation of the conformable double ARA transform and Adomian decomposition methods, to present solutions to the one-dimensional regular and singular conformable fractional coupled Burgers' problem. Additionally, two intriguing examples were given to demonstrate the applicability of the novel approach. Different types of nonlinear time-fractional differential equations with conformable derivatives can be solved using this technique. We want to answer more fractional integral equations and fractional nonlinear problem classes in the future.

*Acknowledgement:*
The authors express their gratitude to the dear referees, who wish to remain anonymous, and the editor for their helpful suggestions.

*References:*

[1] Saadeh, R. A reliable algorithm for solving system of multi-pantograph equations. *WSEAS Trans. Math,* 2022, 21, 792-800.

[2] Gharib, G., & Saadeh, R. Reduction of the self-dual yang-mills equations to sinh-poisson equation and exact solutions. *WSEAS Interact*. 2021, Math, 20, 540-546.

[3] Saadeh, R. Numerical algorithm to solve a coupled system of fractional order using a novel reproducing kernel method. *Alexandria Engineering Journal*, 2021, 60(5), 4583-4591.

[4] J.M. Burgers, A Mathematical Model Illustrating the Theory of Turbulence. Adv. Appl. Mech., 1, 1948, 171–199.

[5] O. Özkan, A. Kurt, On conformable double Laplace transform. *Opt. Quant. Electron., 50*, 2018, 103.

[6] Y. Çenesiz, D. Baleanu, A. Kurt, O.Tasbozan, New exact solutions of Burgers' type equations with conformable derivative. *Wave Random Complex Media*, 27, 2017, 103–116.

[7] J. Liu, G. Hou, Numerical solutions of the space-and time-fractional coupled Burgers equations by generalized differential transform method. *Appl. Math. Comput.*, 217, 2011, 7001–7008.

[8] H. Eltayeb, I. Bachar, A. Kılıçman, On Conformable Double Laplace Transform and One Dimensional Fractional Coupled Burgers' Equation. *Symmetry,* 11, 2019, 417.

[9] J. Nee, J. Duan, Limit set of trajectories of the coupled viscous Burgers' equations. *Appl. Math. Lett.*,11, 1998, 57–61

[10] M. Z. Mohamed, A. E. Hamza, A. K. H. Sedeeg, Conformable double Sumudu transformations an efficient approximation solution to the fractional coupled Burger's equation. *Ain Shams Engineering Journal*, 14(3), 2023, 101879.

[11] R. Saadeh, A. Qazza, A. Burqan, A new integral transform: ARA transform and its properties and applications. *Symmetry*,12, 2020, 925.

[12] R. Saadeh, A. Qazza, A. Burqan, On the Double ARA-Sumudu Transform and Its Applications. *Mathematics*, 10, 2022, 2581.

[13] A. Qazza, A. Burqan, R. Saadeh, R. Khalil, Applications on Double ARA–Sumudu Transform in Solving Fractional Partial Differential Equations. *Symmetry*, 14, 2022, 1817.

[14] E. Salah, A. Qazza, R. Saadeh, A. El-Ajou, A hybrid analytical technique for solving multi-dimensional time-fractional Navier-Stokes system. *AIMS Mathematics*, 8(1), 2023, 1713-1736.

[15] R. Saadeh, Applications of Double ARA Integral Transform. *Computation*, 10, 2022, 216.

[16] R. Saadeh, Application of the ARA Method in Solving Integro-Differential Equations in Two Dimensions. *Computation*, 11, 2023, 4.

[17] A. Qazza, R. Saadeh, E. Salah, Solving fractional partial differential equations via a

new scheme. *AIMS Mathematics*, 8(3), 2023, 5318-5337.

[18] A. Qazza, R. Saadeh, On the analytical solution of fractional SIR epidemic model. *Applied Computational Intelligence and Soft Computing*, 2023, 2023.

[19] R. Rach, On the Adomian (decomposition) method and comparisons with Picard's method. *Journal of Mathematical Analysis and Applications*, 128(2), 1987, 480-483.

[20] M. Younis, A. Zafar, K.U. Haq, M. Rahman, Travelling wave solutions of fractional order coupled Burger's equations by (G0 /G)-expansion method. Am. *J. Comput. Appl. Math.*, 3, 2013, 81.

[21] M.S. Hashemi, Invariant subspaces admitted by fractional differential equations with conformable derivatives. *Chaos Solitons Fractals*, 107, 2018, 161–169.

[22] T. Abdeljawad, On conformable fractional calculus. *J. Comput. Appl. Math.*, 279, 2015, 57–66.

[23] Saadeh, R. (2022). Analytic Computational Method for Solving Fractional Nonlinear Equations in Magneto-Acoustic Waves. *WSEAS Transactions on Fluid Mechanics*, 17, 241-254.

[24] T. Abdeljawad, M. Al-Horani, R. Khalil, Conformable fractional semigroups of operators. *J. Semigroup Theory Appl*. 2015, 2015, 7.

**Contribution of Individual Authors to the Creation of a Scientific Article (Ghostwriting Policy)**

The authors equally contributed to the present research, at all stages from the formulation of the problem to the final findings and solution.

**Sources of Funding for Research Presented in a Scientific Article or Scientific Article Itself**

No funding was received for conducting this study.

**Conflict of Interest**

The authors have no conflicts of interest to declare.